\documentclass[12pt]{amsart}
\usepackage{amssymb}

\newtheorem{theorem}{Theorem}

\newtheorem{prop}[theorem]{Proposition}

\newtheorem{lemma}[theorem]{Lemma}

\newtheorem{cor}[theorem]{Corollary}

\newcommand{\eqa}{\begin{eqnarray}}
\newcommand{\eeqa}{\end{eqnarray}}
\newcommand{\beq}{\begin{equation}}
\newcommand{\eeq}{\end{equation}}
\newcommand{\nn}{\nonumber}

\newcommand{\pal}{\partial}
\newcommand{\cstar}{\buildrel * \over c}

\newcommand{\al}{\alpha}

\newcommand{\pf}{\noindent{\it Proof \ }}

\newcommand{\epf}{$\quad$\hfill
\raisebox{0.11truecm}{\fbox{}}\par\vskip0.4truecm}

\newcommand{\diag}{{\rm diag}\,}
\begin{document}

\title[Duality for Frobenius manifolds]
{On almost duality for Frobenius manifolds}

\author[B. Dubrovin]{Boris Dubrovin}

\address{SISSA\\ Via Beirut 2--4\\ 34014 Trieste\\ Italy
\\ and 
Steklov Math. Institute\\ Moscow}

\email{dubrovin@@sissa.it}

\maketitle

\begin{flushright}
\parbox{9cm}{
\begin{center}
{\it Dedicated
to Sergei Petrovich Novikov \\
on the occasion of his 65th birthday.}
\end{center}
}
\end{flushright}

\begin{abstract}
We present a universal construction of almost duality for Frobenius manifolds.
The analytic setup of this construction is described in details for the case
of semisimple Frobenius manifolds. We illustrate the general considerations by
examples from the singularity theory, mirror symmetry, 
the theory of Coxeter groups and Shephard
groups, from the Seiberg - Witten duality.
\end{abstract}

\section{Introduction}

The beauty of the theory of Frobenius manifold is not only in multiple
connections of it with other branches of mathematics, such as quantum
cohomology, singularity theory, the theory of integrable systems. Even more
amazing is that, some properties 
discovered in the study of particular classes of Frobenius manifolds
often turn out to become {\it universal structures} of the theory thus proving
to be important also for other classes of Frobenius manifolds.

In this paper we describe one of these universal structures. We call it almost
duality; it associates to a given Frobenius manifold a somewhat different
creature we called almost Frobenius manifold. Its crucial part is still
in the WDVV associativity equations; however
the properties of the unity and Euler
vector fields are to be modified.

For the Frobenius structures on the base of universal unfoldings of isolated
hypersurface singularities the duality is based on the wellknown correspondence
between complex oscillatory integrals and periods of closed forms over vanishing
cycles. For quantum cohomologies the duality seems to be closely related to the
mirror construction; at least this is the case in simple examples. In the
setting of the theory of integrable systems the duality generalizes the property
of the classical Miura transformation between KdV and modified KdV equations.

The paper is organized as follows. In Section \ref{crash} we recall necessary
information from the general theory of Frobenius manifolds. In \ref{sec2}
we introduce the almost duality and prove the Reconstruction Theorem that
inverts the duality. The main results of Section \ref{sec3} 
are devoted to the analytic properties of the so-called deformed flat
coordinates on the almost dual to a semisimple Frobenius manifold describing
them in terms of the monodromy data of the latter. Finally in Section \ref{sec4}
we consider examples and applications of almost duality.

{\bf Acknowledgments} This work was 
partially supported by Italian Ministry of Education, Universities and
Researches grant  Cofin2001 ``Geometry
of Integrable Systems''.

\par

\setcounter{equation}{0}
\setcounter{theorem}{0}

\section{A brief introduction into Frobenius manifolds}\label{crash}\par

In this section we will collect necessary definitions and geometric constructions
of the theory of Frobenius
manifolds. We will closely follow \cite{cime, pain}.

\subsection{Deformed flat connection, deformed flat coordinates, and 
spectrum of a Frobenius manifold}

A {\it Frobenius algebra} is a pair $(A, <~,~>)$ where $A$
is a commutative associative algebra with a unity over a field $k$ (we will
consider only the cases $k={\mathbb R},\, {\mathbb C}$) and $<~,~>$ is a $k$-bilinear
symmetric nondegenerate {\it invariant} form on $A$, i.e.,
$$
<x\cdot y, z> = <x, y\cdot z>
$$
for arbitrary vectors $x$, $y$, $z$ in $A$.

{\bf Definition 1.} 
{\it Frobenius structure} of the charge $d$ on the manifold
$M$ is a
structure of a Frobenius algebra on the tangent spaces $T_tM
=(A_t,<~,~>_t)$ depending (smoothly, analytically etc.) on the point $t\in
M$. It must satisfy the following axioms.
 
{\bf FM1.} The metric $<~,~>_t$ on $M$ is flat (but not necessarily
positive definite). Denote $\nabla$ the Levi-Civita connection for the
metric. The unity vector field $e$ must be flat, 
\beq\label{e-flat}
\nabla
e=0.
\eeq
 
{\bf FM2.} Let $c$ be the 3-tensor $c(x,y,z):=<x\cdot y, z>$, $x,\, y,\,
z\in T_tM$. The 4-tensor $(\nabla_w c)(x,y,z)$ must be symmetric in
$x,\, y,\,
z, \, w \in T_tM$.
 
{\bf FM3.} A linear vector field $E\in Vect(M)$ must be fixed on $M$,
i.e. $\nabla\nabla E=0$, such that 
$$
[E, x\cdot y] -[E,x]\cdot y -x\cdot [E,y] = x\cdot y
$$
$$
E<x,y>-<[E,x],y>-<x,[E,y]>=(2-d)<x,y>.
$$

The last condition means that
the derivations
$Q_{Func(M)}:=E, ~~Q_{Vect(M)}:={\rm id}+{\rm ad}_E $
define on the space $Vect(M)$ of vector fields on $M$
a structure of graded Frobenius algebra over the graded ring of functions
$Func(M)$ (see details in \cite{pain}).

Flatness of the metric $<~,~>$ implies local existence of a system of {\it flat
coordinates} $t^1$, \dots, $t^n$ on $M$. We will denote $\eta_{\alpha\beta}$
the constant Gram matrix of the metric  
in these coordinates
$$
\eta_{\alpha\beta}:= \left< {\pal\over \pal t^\alpha}, {\pal\over\pal t^\beta}
\right>.
$$
The inverse matrix $\eta^{\alpha\beta}$ defines the inner product on the
cotangent planes
$$
<dt^\alpha, dt^\beta> =\eta^{\alpha\beta}.
$$
The flat coordinates will be chosen in such a way that the unity $e$ of the
Frobenius algebra coincides with $\pal/\pal t^1$
$$
e={\pal\over\pal t^1}.
$$
In these flat coordinates on $M$ the structure constants of the Frobenius
algebra $A_t=T_tM$ 
\beq\label{wdvv0}
{\pal\over \pal t^\alpha} \cdot {\pal\over\pal t^\beta} =
c_{\alpha\beta}^\gamma(t) {\pal\over\pal t^\gamma}
\eeq
can be locally represented by third derivatives of a function $F(t)$,
$$
c_{\alpha\beta}^\gamma(t) =\eta^{\gamma\epsilon}{\pal^3 F(t)\over \pal
t^\epsilon \pal t^\alpha \pal t^\beta}
$$
also satisfying
\beq\label{wdvv1}
{\pal^3 F(t) \over \pal t^1 \pal t^\alpha \pal t^\beta} \equiv
\eta_{\alpha\beta}.
\eeq
The function $F(t)$ is called {\it potential} of the Frobenius manifold. It is
defined up to adding of an at most quadratic polynomial in $t^1$, \dots, $t^n$.
It satisfies the following system of {\it WDVV associativity equations}
\beq\label{wdvv2}
{\pal^3 F(t)\over \pal t^\alpha \pal t^\beta \pal t^\lambda} \eta^{\lambda\mu}
{\pal^3 F(t)\over \pal t^\mu \pal t^\gamma \pal t^\delta}=
{\pal^3 F(t)\over \pal t^\delta \pal t^\beta \pal t^\lambda} \eta^{\lambda\mu}
{\pal^3 F(t)\over \pal t^\mu \pal t^\gamma \pal t^\alpha}
\eeq
for arbitrary $1\leq\alpha, \, \beta, \, \gamma,\, \delta\leq n$.

The vector field $E$ is called {\it Euler vector field}. In the flat coordinates
it must have the form
$$
E=\left( a^\alpha_\beta t^\beta + b^\alpha\right) {\pal \over \pal t^\alpha}
$$
for some constants $a^\alpha_\beta$, $b^\alpha$ satisying 
$$
a^\alpha_1 =\delta^\alpha_1, ~~b^1=0.
$$
The potential $F(t)$ is a quasihomogeneous function in the following sense
\beq\label{wdvv3}
E\, F= (3-d) F + {1\over 2} A_{\alpha\beta} t^\alpha t^\beta +B_\alpha t^\alpha
+C
\eeq
where $A_{\alpha\beta}$, $B_\alpha$, $C$ are some constants and $d$ is the
charge of the Frobenius manifold.

One of the main geometrical structures of the theory of Frobenius manifolds is
the {\it deformed flat connection}. This is a symmetric affine connection
$\tilde\nabla$
on $M\times {\mathbb C}^*$ defined by the following formulae
\eqa\label{def-con1}
&&
\tilde \nabla_x\, y = \nabla_x y + z\, x\cdot y, ~~x, y\in TM, ~~z\in{\mathbb
C}^* \\
&&
\tilde \nabla _{d\over dz} y = \pal_z y +E\cdot y -{1\over z} {\mathcal V} y
\nn\\
&&
\tilde\nabla_x {d\over dz} = \tilde\nabla_{d\over dz} {d\over dz} =0
\label{def-con2}
\eeqa
where
\beq\label{2-7-4b}
{\mathcal V} := {2-d\over 2} -\nabla E
\eeq
is an antisymmetric operator on the tangent bundle $TM$ w.r.t. $<\ ,\ >$,
$$
<{\mathcal V}\, x, y> =-<x, {\mathcal V}\, y>.
$$
Observe that the unity vector field $e$ is an
eigenvector of this operator with the eigenvalue
$$
{\mathcal V} e = -{d\over 2} e.
$$
Vanishing of the curvature of the connection $\tilde \nabla$ is essentially
equivalent to the axioms of Frobenius manifold.

{\bf Definition 2.} A function $f=f(t;z)$ defined on an open subset in $M\times
{\mathbb C}^*$ is called $\tilde\nabla$-{\it flat} if
\beq\label{f-flat}
\tilde \nabla df =0.
\eeq

Introducing the row vector $\xi=(\xi_1, \dots, \xi_n)$ where
$$
\xi_\al: = {\pal f\over \pal t^\al}
$$
we represent the equations for $\tilde\nabla$-flatness in the form
\eqa\label{xi-flat1}
&&
\pal_\al \xi = z\, \xi \, C_\al(t)
\\
&&
\pal_z \xi = \xi \, \left( {\mathcal U}(t) -{{\mathcal V}\over z} \right).
\label{xi-flat2}
\eeqa
Here we introduce the
matrix ${\mathcal U}(t)$ of multiplication by the Euler vector field
\beq\label{calu}
{\mathcal U}_\beta^\alpha (t) := E^\epsilon(t) c_{\epsilon \beta}^\alpha (t)
\eeq
and the matrix $C_\al(t)$ of multiplication by $\partial_\alpha:= \pal / \pal
t^\al$
\beq\label{big-c}
\left( C_\alpha(t)\right)_\gamma^\beta : = c_{\alpha\gamma}^\beta(t).
\eeq
The same system can be rewritten for the column vector $y=(y^1, \dots, y^n)^T$,
$$
y^\al = \eta^{\al\beta}\xi_\beta
$$
in the form
\eqa\label{lin-system1}
&&
\pal_\al y = z\, C_\al (t) y
\\
&&
\pal_z y = \left( {\mathcal U}(t) +{{\mathcal V}\over z} \right)\,y.
\label{lin-system2}
\eeqa

Because of vanishing of the torsion and curvature of the connection
$\tilde\nabla$ there locally exist, on $M\times
{\mathbb C}^*$ $n$ independent flat functions $\tilde t^1(t;z)$, \dots, 
$\tilde t^n(t;z)$. They are called {\it deformed flat coordinates} on a
Frobenius manifold. The analytic
properties of deformed flat coordinates as multivalued functions of 
$z\in {\mathbb C}^*$ can be used for describing moduli of semisimple 
Frobenius manifolds (see the next Section). Here we briefly describe 
the behaviour near the regular singularity $z=0$ of a particular basis
of deformed flat coordinates (see details in \cite{pain}).

To fix a system of the deformed flat coordinates we are to choose a basis in the space
of solutions to the system (\ref{lin-system1}), (\ref{lin-system2}). 
Such a basis corresponds 
to a choice of a representative in the equivalence class of {\it normal forms} of the system 
(\ref{lin-system2}) near $z=0$ (see details in \cite{pain}). 
The parameters of such a normal form are called {\it spectrum} 
of the Frobenius manifold. Let us first recall the description of the 
parameters.

{\bf Definition 3.} The {\it spectrum} of a Frobenius manifold is a quadruple
$\left( V,<~,~>, \hat\mu, R\right)$ where
$V$ is a $n$-dimensional linear space over ${\mathbb C}$ equipped 
with a symmetric non-degenerate bilinear form $<\ ,\ >$, 
semisimple antisymmetric 
linear operator $\hat \mu: V\to V$, $<\hat\mu\, a, b>+<a, \hat\mu\,b>=0$
and a nilpotent linear operator $R: V\to V$ satisfying the following properties.
First, 
\beq\label{R2}
R^* =-e^{\pi\,i\,\hat\mu} R e^{-\pi\,i\,\hat\mu}
\eeq
Observe the following consequence of (\ref{R2})
\beq\label{levelt51}
R e^{2\pi\,i\,\hat\mu} = e^{2\pi\,i\,\hat\mu} R.
\eeq
In particular, $R$ leaves invariant the eigensubspaces of the operator
$e^{2\pi\,i\,\hat\mu}$. It {\it does not} leave invariant the eigensubspaces
$V_\mu$ of the operator $\hat\mu$. However,
\beq\label{dec}
R\, V_\mu \subset \oplus_{m\in {\mathbb Z}} V_{\mu +m}
\eeq
The crucial condition in the definition of the spectrum is that,
the operator $R$ must also be {\it $\hat\mu$-nilpotent}, i.e., in the
decomposition (\ref{dec}) only nonnegative integers $m$ are present.
We define the components of the operator $R$
\eqa\label{comp-r}
&&
R=R_0+R_1+R_2+\dots
\nn\\
&&
R_m V_\mu \subset V_{\mu+m}, ~~{\rm for} ~{\rm any}~ \mu \in {\rm Spec}\,
\hat\mu.
\eeqa

By the construction the operator $R$ satisfies
\beq\label{R3}
z^{\hat\mu} R z^{-\hat\mu} = R_0 + R_1 z + R_2 z^2+\dots.
\eeq
Observe also the following useful identity
\beq\label{R6}
[\hat\mu, R_k] = k\, R_k, ~~k=0, 1, \dots .
\eeq
Any polynomial of the matrices $R_k$ can be uniquely decomposed as
follows
\beq\label{comp}
P(R_0, R_1, \dots)= [P(R_0, R_1, \dots)]_0+[P(R_0, R_1, \dots)]_1+\dots
\eeq
\beq\label{comp-b}
z^{\hat\mu} [ P(R_0, R_1, \dots)]_m z^{-\hat\mu} = z^m [P(R_0, R_1,
\dots)]_m.
\eeq 
The last restriction for the spectrum is that, the eigenvector $e\in V$
of $\hat\mu$ must satisfy $R_0\,e=0$.

We will now explain how to associate a 5-tuple $(V, <\ ,\ >, \hat\mu, R, e)$ 
to a Frobenius manifold. The linear space $V$ with a symmetric nondegenerate bilinear form $<\ ,\ >$ 
and a vector $e\in V$ have already been constructed above.
Denote $\hat\mu:V\to V$ the semisimple part of the operator ${\mathcal V}$, i.e., 
\beq\label{mu}
\hat\mu:= \oplus_{\mu\in Spec\, {\mathcal V}} \mu\, P_\mu
\eeq
where $P_\mu:V\to V_\mu$ is the projector of $V$ onto the root subspace of 
${\mathcal V}$
$$
V=\oplus_{\mu\in Spec \, {\mathcal V}} V_\mu, ~~~V_\mu := {\rm Ker}\, ({\mathcal V}-\mu\cdot 1)^n,
$$
$P_\mu (V_{\mu'})=0$ for $\mu\neq \mu'$, $P_\mu |_{V_\mu}={\rm id}_{V_\mu}$.
Clearly the operator $\hat\mu$ is antisymmetric, $\hat\mu^*=-\hat\mu$. Denote $R_0$ the nilpotent part of ${\mathcal V}$
$$
{\mathcal V}=\hat\mu+R_0.
$$ 
Other operators $R_1$, $R_2$, \dots are not determined by ${\mathcal V}$ only. They appear only
in
presence of resonances, i.e., pairs of eigenvalues $\mu$, $\mu'$ of ${\mathcal V}$ 
such that $\mu-\mu' \in {\mathbb Z}_{>0}$ (see details in \cite{pain}).

Let us choose a basis $e_1, \dots, e_n$ in $V$ such that $e_1=e$. The matrices
of the linear operators $\hat\mu$ and $R$ we will denote by the same symbols.

\begin{theorem} For a sufficiently small ball $B\in M$ there exists a
fundamental matrix of solutions to the system (\ref{lin-system1}), 
(\ref{lin-system2}) of the form
\beq\label{theta}
y(t;z) = \Theta(t;z)z^{\hat\mu}z^R=\sum_{p\geq 0} \Theta_p(t)z^{p+\hat\mu} z^R
\eeq
such that
the matrix valued function $\Theta(t;z):V\to V$ is analytic on $B\times {\mathbb C}$ satisfying
\beq\label{normalize-theta}
\Theta(t;0)\equiv 1
\eeq
\beq\label{orthogonal}
\Theta^*(t;-z)\Theta(t;z)\equiv 1.
\eeq
\end{theorem}

By the construction the columns of the fundamental matrix $Y(t;z)$ are gradients
of a certain system of deformed flat coordinates
\beq\label{vtilde}
(\tilde t(t; z), \dots, \tilde t(t;z)) =(\theta_1(t;z),
\dots,\theta_n(t;z))z^{\hat\mu} z^R
\eeq 
where
\beq\label{theta-alpha-p}
\theta_\alpha(t;z) =\sum_{p=0}^\infty \theta_{\alpha,p}(v) z^p, ~~\alpha=1,
\dots, n.
\eeq

{\bf Definition 4.} We will call (\ref{vtilde}) {\it Levelt system of deformed
flat coordinates} on $M$ at $z=0$.

The ambiguity in the choice of the Levelt basis of deformed flat coordinates
is described in \cite{pain}.

\subsection{Canonical coordinates on semisimple Frobenius manifolds}
\label{sec-3-10-4a}\par

{\bf Definition 5.} The Frobenius manifold $M$ is called {\it semisimple} 
if the family of
$n$-dimensional algebras (\ref{wdvv0}) is
semisimple for any $t=(t^1, \dots, t^n)\in M_{s\, s}$ for an open dense subset
$M_{s\, s}\subset M$.

In this Section we summarize, following \cite{npb, cime} 
a very efficient technique of working with
semisimple Frobenius manifolds based on introduction of {\it canonical
coordinates}.

Let $M$ be a semisimple Frobenius manifold. Denote $M_{s\, s}\subset M$
the open dense subset in $M$ consisting of all points $t\in M$ s.t.
the operator ${\mathcal U}(t)$ of multiplication by the Euler vector field
$$
{\mathcal U}(t)=E(t)\cdot \, : T_t M\to T_t M
$$
has simple spectrum. (Actually, the subset $M_{s\, s}\subset M$ could be
slightly smaller than the set of all points of semisimplicity of the algebra
(\ref{wdvv0}). It can be shown however that, for an analytic Frobenius manifold
this is still an open dense subset.) Denote $u_1(t)$, \dots, $u_n(t)$ the eigenvalues
of this operator, $t\in M_{s\, s}$. 

\begin{theorem} The mapping
$$
M_{s\, s} \to \left({\mathbb C}^n\setminus \cup_{i<j}(u_i=u_j) \right)
/S_n, ~~t\mapsto (u_1(t), \dots, u_n(t))
$$
is an unramified covering. 
Therefore one can use the eigenvalues as local
coordinates on $M_{s\, s}$. In these coordinates the multiplication table of the
algebra (\ref{wdvv0}) becomes 
\beq\label{idem}
{\pal \over \pal u_i} \cdot {\pal \over \pal u_j} = 
\delta_{ij} {\pal \over \pal u_i}.
\eeq
The basic idempotents $\pal/ \pal u_i$ are pairwise orthogonal
$$
\left< {\pal \over \pal u_i} , {\pal \over \pal u_j}\right> =0, ~i\neq j.
$$
\end{theorem}

Observe that we violate
the indices convention labelling the canonical coordinates by subscripts.
We will {\it never} use summation over repeated indices when working in the canonical
coordinates.

Choosing locally branches of the square roots
\beq\label{psi-1}
\psi_{i1}(u):= \sqrt{<\pal/\pal u_i, \pal/\pal u_i>}, ~i=1, \dots, n
\eeq
we obtain a transition matrix $\Psi =(\psi_{i\alpha}(u))$
from the basis $\pal /\pal t^\alpha$ to the
orthonormal basis
\beq\label{basis}
f_1=\psi_{11}^{-1}(u) {\pal \over \pal u_1} ,\, f_2=\psi_{21}^{-1}(u) {\pal \over \pal u_2} ,
\dots ,\, f_n= \psi_{n1}^{-1}(u) {\pal \over \pal u_n} 
\eeq
of the normalized idempotents 
\beq\label{mat-psi}
{\pal\over \pal t^\alpha} =\sum_{i=1}^n {\psi_{i\alpha}(u)\over
\psi_{i1}(u)}{\pal\over \pal u_i}.
\eeq
Equivalently, the Jacobi matrix has the form
\beq
{\pal u_i\over \pal t^\alpha} ={\psi_{i\alpha}\over\psi_{i1}}.
\eeq
The matrix $\Psi(u)$ satisfies orthogonality condition
\beq\label{orth}
\Psi^*(u) \Psi(u) \equiv 1, ~~\Psi^*:= \eta^{-1} \Psi^T \eta, ~~\eta = (\eta_{\alpha\beta}),
~~\eta_{\alpha\beta}:= \left< {\pal \over \pal t^\alpha}, {\pal \over \pal
t^\beta}\right>.
\eeq
In this formula $\Psi^T$ stands for the transposed matrix.
The lengths (\ref{psi-1}) coincide with the first column of this matrix. So,
the metric $<\ ,\ >$ in the canonical coordinates reads
\beq\label{metric-d}
<\ ,\ >=\sum_{i=1}^n \psi_{i1}^2(u) du_i^2.
\eeq
The inverse Jacobi matrix can be computed using (\ref{orth}). This gives
\beq\label{inv-jac}
{\pal t_\alpha\over \pal u_i} = \psi_{i1} \psi_{i\al}, \quad t_\al:=
\eta_{\al\beta} t^\beta.
\eeq

Denote $V(u)=\left( V_{ij}(u)\right)$ the matrix of the antisymmetric operator 
${\mathcal V}$ 
(\ref{2-7-4b}) w.r.t. the orthonormal frame (\ref{basis})
\beq\label{mat-v}
V(u) := \Psi(u)\, {\mathcal V}\, \Psi^{-1}(u).
\eeq

The matrix $V(u)$
satisfies the following system of commuting time-dependent
Hamiltonian flows on the Lie algebra $so(n)$ equipped with the standard Lie -
Poisson bracket 
\beq\label{eq-per-v}
{\pal V\over \pal u_i} =  \{ V, H_i(V; u)\}, ~~i=1, \dots, n
\eeq
with quadratic Hamiltonians
\beq\label{iso-ham}
H_i(V; u) ={1\over 2} \sum_{j\neq i} {V_{ij}^2\over u_i -u_j}.
\eeq
The matrix $\Psi(u)$ satisfies
\beq\label{eq-psi}
{\pal\Psi\over \pal u_i} = V_i(u) \Psi, ~~ V_i(u):= {\rm ad}_{E_i} {\rm ad}_U
^{-1} (V(u)).
\eeq
Here 
\beq\label{mat-u}
U={\rm diag}\, (u_1, \dots, u_n)=\Psi\, {\mathcal U} \Psi^{-1},
\eeq 
the matrix unity $E_i$ has the entries
\beq\label{entr}
(E_i)_{ab}= \delta_{ai} \delta_{ib}.
\eeq


The system (\ref{eq-per-v}) coincides with the equations of isomonodromy
deformations of the following linear differential operator with rational
coefficients
\beq\label{isomono}
{dY\over  dz} = \left( U+{V\over z}\right) \, Y.
\eeq
The latter is nothing but the last
component of the deformed flat connection (\ref{def-con2}) written in the
orthonormal frame (\ref{basis}). Namely, the solutions to (\ref{isomono}) are
related to the gradients of the $\tilde\nabla$-flat functions satisfying 
(\ref{lin-system2})  by
\beq\label{kalibr} 
Y=\Psi^{-1} y.
\eeq
The first part (\ref{lin-system1}) after the gauge transformation (\ref{kalibr})
reads
\beq\label{isomono1}
\pal_i Y = (z\, E_i +V_i) Y.
\eeq
In particular the fundamental matrix (\ref{theta}) gives the fundamental matrix
of solutions to (\ref{isomono}), (\ref{isomono1}) of the form
\eqa\label{theta1}
&&
Y_0 (u;z)= \Psi^{-1}(u) \Theta(t(u);z)z^{\hat\mu}z^R
\nn\\
&&
=\sum_{p\geq 0} \Theta_p(t(u))z^{p+\hat\mu} z^R.
\eeqa

The integration of (\ref{eq-per-v}),
(\ref{eq-psi}) and, more generally, the
reconstruction of the Frobenius structure can be reduced to a solution
of certain Riemann - Hilbert problem \cite{pain}.

\subsection{Intersection form, discriminant, and periods
of a Frobenius manifold.}\label{int}\par

We now recall some important points of the theory of intersection
form on a Frobenius manifold and of the corresponding period mapping
(see \cite{cime, tan, pain}). Intersection form is a symmetric  bilinear form on 
$T^*M$. It is defined by the formula
\beq\label{int-form}
\left( \omega_1, \omega_2\right) =i_E (\omega_1\cdot \omega_2),
~~\omega_1, \, \omega_2 \in T_t^* M.
\eeq
The product $\omega_1\cdot \omega_2$ of 1-forms is induced from the
product of vectors in $T_tM$ by the isomorphism
$$
<~,~>: T_tM\to T_t^*M.
$$
This means that, in the flat coordinates for $<~,~>$,
\beq\label{prod}
dt^\alpha\cdot dt^\beta =c^{\alpha\beta}_\gamma(t) dt^\gamma, ~~
c^{\alpha\beta}_\gamma(t) = \eta^{\alpha\lambda}\eta^{\beta{\mathcal V}}
{\partial^3 F(t)\over \partial t^\lambda \partial t^{\mathcal V} \partial
t^\gamma}=\eta^{\al\lambda} c_{\lambda\gamma}^\beta(t).
\eeq
For the Gram matrix 
$$
g^{\alpha\beta}(t) := \left( dt^\alpha, dt^\beta\right)
$$
of the bilinear form (\ref{int-form}) one obtains
\eqa\label{int-form1}
&&
g^{\alpha\beta}(t) = E^\epsilon(t) c_\epsilon^{\alpha\beta}(t) 
\nn\\
&& 
= F^{\alpha\beta}(t) -{\mathcal V}_\rho^\alpha F^{\rho\beta}(t) -F^{\alpha\rho}(t)
{\mathcal V}_\rho^\beta +A^{\alpha\beta}
\eeqa
Here
$$
F^{\alpha\beta}(t) =\eta^{\alpha\alpha'} \eta^{\beta\beta'}{\partial^2
F(t)\over \partial t^{\alpha'} \partial t^{\beta'}},
$$
$$
A^{\alpha\beta}=\eta^{\alpha\alpha'} \eta^{\beta\beta'} A_{\alpha'\beta'}
$$
and the constant matrix $\left(A_{\alpha\beta}\right)$ was defined in
(\ref{wdvv3}).

\begin{prop}\label{prop1-1} For any complex parameter $\lambda$ denote
$\Sigma_\lambda \subset M$ the subset
\beq\label{discr}
\Sigma_\lambda :=\left\{ t\in M ~ | ~\det \left[ (~,~)_t -\lambda <~,~>_t
\right] =0 \right\} .
\eeq
It is a proper analytic subset (i.e., $\det \left( g^{\alpha\beta}(t)
-\lambda \eta^{\alpha\beta}\right)$ does not vanish identically on $M$).
On $M\setminus \Sigma_\lambda$ the inverse matrix
$$
g_{\alpha\beta}(t;\lambda): =  \left( g^{\alpha\beta}(t)
-\lambda \eta^{\alpha\beta}\right)^{-1}
$$
defines a flat metric that we denote $(~,~)_\lambda$. The Levi-Civita
connection of this metric has the Christoffel coefficients
$\Gamma_{\beta\gamma}^\alpha(t;\lambda)$ of the form
\beq\label{christ0}
\Gamma_{\beta\gamma}^\alpha(t;\lambda) =- g_{\beta\rho}(t;\lambda)
\Gamma_\gamma^{\rho\alpha}(t)
\eeq
where
\beq\label{christ}
\Gamma_\gamma^{\alpha\beta}(t) = c_\gamma^{\alpha\epsilon}(t) \left(
{1\over 2} -{\mathcal V}\right) _\epsilon^\beta.
\eeq
\end{prop}

The functions $\Gamma_\gamma^{\alpha\beta}(t)$ are defined everywhere on
$M$ (but not only on $M\setminus \Sigma_\lambda$). We will call them {\it
contravariant Christoffel coefficients} of the metric $(~,~)_\lambda$.

{\bf Remark 1.} The formulae of Proposition imply that the metrics
$(~,~)$ and $<~,~>$ on $T^*M$ form {\it flat pencil} \cite{cime, tan}. This means that:

1). The contravariant Christoffel coefficients for an arbitary linear
combination $(~,~)-\lambda <~,~>$ on $T^*M$ are 
\beq\label{pen1}
{\Gamma_.^{\cdot \,
\cdot}}_{(~,~)} -\lambda {\Gamma_.^{\cdot \,
\cdot}}_{<~,~>}
\eeq
where ${\Gamma_.^{\cdot \,
\cdot}}_{(~,~)}$ and ${\Gamma_.^{\cdot \,
\cdot}}_{<~,~>}$ are the contravariant Christoffel coefficients 
of the Levi-Civita connections $\nabla^{(~,~)}$ and $\nabla^{<~,~>}$
for the metrics $(~,~)$ and $<~,~>$ resp.

2). The linear combination $(~,~)-\lambda<~,~>$ is a flat metric on $T^*_tM$,
$t\in M\setminus \Sigma_\lambda$ for any
$\lambda\in {\mathbb C}$.

3). The flat pencil is said to be {\it quasihomogeneous of the charge $d$}
if a function $f$ exists such that the Lie derivatives of the metrics 
along the vector
fields
\beq\label{pen2}
 E:= \nabla^{(~,~)} f, ~~e:= \nabla^{<~,~>}f
\eeq
have the form
\eqa\label{pen3}
&&
Lie_E (~,~) =(d-1) (~,~), ~~ Lie_e (~,~) =<~,~>,
\\
&&
Lie_E <~,~> =(d-2) <~,~>, ~~ Lie_e <~,~> =0
\label{pen4}
\eeqa
for some constant $d\in {\mathbb C}$
and the commutator of the vector fields equals
\beq\label{pen5}
[e, E] = e.
\eeq
In the case of Frobenius manifolds the connection $\nabla^{<~,~>}$ coincides
with $\nabla$, the connection $\nabla^{(~,~)}$ defined by (\ref{christ}) for
$\lambda=0$
will be denoted $\nabla_*$.
The function $f$ has the form
$$
f=t_1 \equiv \eta_{1\alpha}t^\alpha.
$$
As it was shown in \cite{cox} (see also \cite{cime, tan}) existence of a flat quasihomogeneous
pencil with certain restrictions for the eigenvalues of the operator
$\nabla^{<~,~>} E$ can be
used for an alternative axiomatization of Frobenius manifolds.

{\bf Definition 6.} A function $p=p(t;\lambda)$ is called {\it $\lambda$-period} of the
Frobenius manifold if it satisfies
\beq\label{period}
(\nabla_* - \lambda \nabla) dp =0.
\eeq

Due to Proposition \ref{prop1-1}, on the universal covering of
$$
M\times {\mathbb C} \setminus \cup_\lambda \lambda \times \Sigma_\lambda
$$
there exist $n$ independent $\lambda$-periods $p^1(t;\lambda)$, \dots, $p^n(t;\lambda)$.
They give a system of flat coordinates for the metric $(~,~)-\lambda<~,~>$.
More precisely,

\begin{cor} The flat coordinates $p^1(t;\lambda)$, \dots,
$p^n(t;\lambda)$ of the metric $(~,~)-\lambda <~,~>$ on a sufficiently
small domain in $M\setminus \Sigma_\lambda$ can be determined from the
system of linear differential equations
\eqa\label{g-m}
&&
\left( g^{\alpha\epsilon}(t) -\lambda\, \eta^{\alpha\epsilon}\right)
\partial_\beta \xi_\epsilon + c_\beta^{\alpha\rho}(t) \left( {1\over
2}-{\mathcal V}\right) _\rho^\epsilon \xi_\epsilon =0, ~\alpha, \, \beta = 1,
\dots, n
\nn\\
&&
\xi_\epsilon ={\partial x(t; \lambda)\over \partial t^\epsilon},
~~\epsilon = 1, \dots, n.
\eeqa
A full system of independent coordinates $p^a=p^a(t;\lambda)$ is obtained from a fundamental
system of solutions $\xi^a(t;\lambda)=\left( \xi^a_\epsilon(t;\lambda)\right)$,
 $a=1, \dots, n$, of the
linear system (\ref{g-m}). In these coordinates the Gram matrix of the bilinear
form
(\ref{int-form}) is a constant one
\beq\label{int-form2}
G^{ab} \equiv \left( dp^a, dp^b\right)-\lambda <dp^a, dp^b> 
=\left[ g^{\alpha\beta}(t)
-\lambda \, \eta^{\alpha\beta}\right]
\xi^a_\alpha(t;\lambda)\xi^b_\beta(\lambda,t)
\eeq
The dependence of the new flat coordinates on $t$, $\lambda$ can be chosen
in such a way that the partial derivatives $\xi_\epsilon(t;\lambda)$
satisfy also a system of linear differential equations with rational
coefficients
\beq\label{g-m1}
\left( g^{\alpha\epsilon}(t) -\lambda\, \eta^{\alpha\epsilon}\right)   
{\partial\xi_\epsilon\over \partial \lambda} -
\eta^{\alpha\rho}\left( {1\over
2}-{\mathcal V}\right) _\rho^\epsilon \xi_\epsilon =0, ~\alpha = 1,\dots, n.
\eeq
The corresponding flat coordinates have the form 
\beq\label{l-shift}
p^a(t;\lambda) =\hat p^a\left(t^1-\lambda, t^2, \dots, t^n\right)   
\eeq
where $\hat p^a(t)$ are flat coordinates of the metric $(~,~)$, $t\in M\setminus
\Sigma_0$.
\end{cor}

We will omit the hat over $p^a(t)$ in the notations for the flat
coordinates of the intersection form.

Let us rewrite the system (\ref{g-m}), (\ref{g-m1}) in matrix notations. 
The equations (\ref{g-m}), (\ref{g-m1}) for the row-vector $\xi =
\left(\partial_1 p,
\partial_2p, \dots, \partial_np\right)$ read
\eqa\label{g-m0}
&&
\partial_\alpha \xi \cdot ({\mathcal U}-\lambda) =\xi \left( {\mathcal V}-{1\over 2}
\right) C_\alpha,
\\
&&
\partial_\lambda \xi \cdot ({\mathcal U}-\lambda) =\xi \left({1\over
2}-{\mathcal V}\right) .
\label{g-m01}
\eeqa
The matrices ${\mathcal U}(t)$, $C_\al(t)$ were defined in (\ref{calu}),
(\ref{big-c}).
Observe that, for $d\neq 1$, one can reconstruct the $\lambda$-period
$p(t;\lambda)$ knowing its gradient $\xi(t;\lambda)$ using

\begin{lemma}\label{lemma1-1}  Let $\xi=\left( \xi_\alpha(t;\lambda)\right)$ be
an arbitrary solution of the system (\ref{g-m0}), (\ref{g-m01}), and $d\neq 1$. Then the
function
\beq\label{reco}
p(t;\lambda) ={2\over 1-d} i_{E-\lambda\,e} \xi \equiv
{2\over 1-d} \left[ E^\epsilon (t) \xi_\epsilon (t;\lambda) -\lambda\,
\xi_1(t;\lambda)\right]
\eeq
satisfies (\ref{g-m}).
\end{lemma}

\pf Multiplying (\ref{g-m0}) by $E^\alpha(t) -\lambda\, \delta_1^\alpha$ and
taking the sum over $\alpha$ one obtains
$$
\left( E^\alpha(t) -\lambda\,\delta^\alpha_1 \right) \partial_\alpha \xi
=\xi \left( {\mathcal V} -{1\over 2}\right).
$$
Using
$$
\nabla E={1-d\over 2} + {1\over 2} -{\mathcal V}
$$
and closedness of the 1-form $\xi$ we rewrite the above equation in the
form
$$
d\, \left( i_{E-\lambda\,e} dp\right)={1-d\over 2} dp.
$$
So
$$
i_{E-\lambda\, e} dp ={1-d\over 2}p +{\rm const}.
$$
Doing a shift along $p$ we kill the constant if $d\neq 1$. Lemma is
proved.

{\bf Definition 7.} The subset 
$$
\Sigma:= \Sigma_0\subset M
$$
is called {\it discriminant} of the Frobenius manifold $M$. Any function
$p=p(t)$ on $M\setminus \Sigma$ satisfying (\ref{g-m0}), (\ref{g-m01}) 
with $\lambda=0$ is called {\it period} of
the Frobenius manifold. A system of $n$ independent periods $p^1(t)$, \dots,
$p^n(t)$ gives
flat
coordinates of the intersection form. They determine a local isometry
of the complement $M\setminus\Sigma$ to the complex Euclidean space ${\mathbb
C}^n$ equipped with the quadratic form $G=\left(G^{ab}\right) = \left( (dp^a,
dp^b)\right)$. This map
is called {\it  period mapping} of the Frobenius manifold.

Sometimes we will call it {\it even} period mapping to distinguish from the odd
one to be introduced in Section 2.

Analytic continuation of the period mapping ${\bf p}(t):=\left( p^1(t), \dots,
p^n(t)\right)$ is a single-valued analytic vector-function on the
universal covering of $M\setminus\Sigma$. Continuing this function along
a closed loop $\gamma$ on $M\setminus\Sigma$ one obtains a new system of
flat coordinates
\beq\label{mono}
{\bf p}(t) \mapsto {\bf p}(t)\, M_\gamma +a_\gamma
\eeq
where $M_\gamma$ is a $n\times n$ matrix satisfying
$$
M_\gamma G\,M_\gamma^T = G
$$
and $a_\gamma$ is a constant vector. The matrix $M_\gamma$ and the vector
$a_\gamma$ depend only on the homotopy class of the loop $[\gamma]\in
\pi_1(M\setminus\Sigma)$. For $d\neq 1$, due to the quasihomogeneity
(\ref{reco})
of the components of the period mapping, one may assume that $a_\gamma=0$.
We will mainly consider here only the case $d\neq 1$.

The following simple statement will be useful later on.

\begin{lemma}\label{emma} If, for $d\neq 1$, the flat coordinates $p^1(t)$, \dots,
$p^n(t)$ of the intersection form are chosen in such a way that
\beq\label{q-h}
Lie_Ep={1-d\over 2} p
\eeq
and 
$$
(dp^a, dp^b) =G^{ab}, ~~(G_{ab})=(G^{ab})^{-1},
$$
then
\beq\label{t1}
t_1:= \eta_{1\alpha}t^\alpha ={1-d\over 4} G_{ab}p^a p^b.
\eeq
\end{lemma}

\pf From (\ref{int-form}) we have the following expression for the gradient of
the function $t_1$ w.r.t. the metric $(~,~)$
\beq\label{q-t1}
\nabla_*t_1 =E.
\eeq
Equation (\ref{q-h}) implies that
\beq\label{q-e}
E={1-d\over 2} p^a {\partial\over \partial p^a}.
\eeq
Rewriting (\ref{q-t1}),  in the flat coordinates $p^1$, \dots, $p^n$ we obtain
$$
G^{ab} {\partial t_1\over\partial p^b} ={1-d\over 2} p^a.
$$
This proves Lemma.

The representation (for $d\neq 1$)
\beq\label{mono-rep}
\pi_1(M\setminus\Sigma) \to O({\mathbb C}^n, G), ~~\gamma \mapsto
M_\gamma^{-1}
\eeq
is called {\it monodromy representation} of the Frobenius manifold.

\setcounter{equation}{0}
\setcounter{theorem}{0}

\section{Dual (almost) Frobenius manifold and its deformed flat
connection}\label{sec2}

Here we will construct a transformation
associating to an arbitrary Frobenius manifold $M$ a new structure on
$M^*=M\setminus\Sigma$ satisfying all the axioms of a Frobenius manifold but
constancy of the unity (cf. \cite{cime}, Remark 4.2). 

For $t\in M\setminus\Sigma$ we define a new multiplication of tangent
vectors $u, \, v\in T_tM$ by
\beq\label{star}
u*v :={u\cdot v\over E}.
\eeq

\begin{prop}\label{prop1-2} The multiplication (\ref{star}) together with the
intersection form $(~,~)$, the unity = the Euler vector field = $E$
satisfies all the axioms of Frobenius manifold but (\ref{e-flat}).
\end{prop}

\pf Associativity of the multiplication is obvious. As it follows from
the
definition, the superposition
$$
f: T_tM {\buildrel {(~,~)^{-1}} \over \longrightarrow} T_t^* M 
{\buildrel {<~,~>}\over \longrightarrow} T_t M
$$
(here $(~,~)$, $<~,~>$ are considered as bilinear forms on $T_t^*M$)
establishes an isomorphism of the algebras
$$
f(u*v) = f(u)\cdot f(v).
$$
In other words, on the cotangent planes both the multiplications are given
by the same formula (\ref{prod}). From this we derive invariance of the
intersection form w.r.t. the new multiplication, i.e., the symmetry of
the expression
$$
g^{\alpha\beta}_\epsilon g^{\epsilon\gamma} =i_E \left( dt^\alpha\cdot
dt^\beta\cdot dt^\gamma \right)
$$
w.r.t. $\alpha$, $\beta$, $\gamma$.

Next, we are to prove the symmetry of the covariant derivatives
$$
\nabla_*^\gamma c^{\alpha\beta}_\rho =g^{\gamma\epsilon}
\partial_\epsilon c^{\alpha\beta}_\rho
-\Gamma_\epsilon^{\gamma\alpha} c_\rho^{\epsilon\beta}
-\Gamma_\epsilon^{\gamma\beta}c_\rho^{\alpha\epsilon}
+\Gamma_\rho^{\gamma\epsilon} c_\epsilon^{\alpha\beta}
$$
w.r.t. $\alpha$, $\beta$, $\gamma$. Here $\nabla_*$ is the Levi-Civita
connection for the intersection form. Rewriting the first term in the
r.h.s. as
$$
g^{\gamma\epsilon}\partial_\epsilon c^{\alpha\beta}_\rho = 
g^{\gamma\epsilon}\partial_\rho c^{\alpha\beta}_\epsilon=
\partial_\rho \left( 
c^{\alpha\beta}_\epsilon g^{\epsilon\gamma}\right)
-c_\epsilon^{\alpha\beta}\, \left(
\Gamma_\rho^{\gamma\epsilon}+\Gamma_\rho^{\epsilon\gamma}\right)
$$
(here we use the condition $\nabla_* g^{\alpha\beta}=0$), we obtain
$$
\nabla_*^\gamma c_\rho^{\alpha\beta} =\partial_\rho \left(
c^{\alpha\beta}_\epsilon g^{\epsilon\gamma}\right)
-\Gamma_\epsilon^{\gamma\alpha} c_\rho^{\epsilon\beta}
-\Gamma_\epsilon^{\gamma\beta}c_\rho^{\alpha\epsilon}
-c_\epsilon^{\alpha\beta} \Gamma_\rho^{\epsilon\gamma}
$$
$$
=\partial_\rho \left(
c^{\alpha\beta}_\epsilon
g^{\epsilon\gamma}\right)-c_\epsilon^{\gamma\lambda}c_\rho^{\epsilon\beta}
\left({1\over 2} -{\mathcal V}\right) _\lambda^\alpha
-c_\epsilon^{\gamma\lambda}c_\rho^{\alpha\epsilon} \left({1\over 2}
-{\mathcal V}\right)_\lambda^\beta -c_\epsilon^{\alpha\beta}
c_\rho^{\epsilon\lambda} \left({1\over 2} -{\mathcal V}\right)_\lambda^\gamma.
$$
Using associativity we recast the last expression into the form
$$
= \partial_\rho \left(
c^{\alpha\beta}_\epsilon
g^{\epsilon\gamma}\right)-\left[
c_\epsilon^{\gamma\beta}c_\rho^{\epsilon\lambda}
\left({1\over 2} -{\mathcal V}\right) _\lambda^\alpha
+c_\epsilon^{\alpha\gamma}c_\rho^{\epsilon\lambda} \left({1\over 2}
-{\mathcal V}\right)_\lambda^\beta +c_\epsilon^{\alpha\beta}
c_\rho^{\epsilon\lambda} \left({1\over 2} -{\mathcal V}\right)_\lambda^\gamma
\right].
$$
Using the symmetry of the first term we clearly see that the expression is
symmetric in $\alpha$, $\beta$, $\gamma$.

It remains to prove that $E$ is the Euler vector field also for the new
algebra structure, and that it also plays the role of the unity for
the new multiplcation. The first statement follows from the formulae
for the Lie derivatives
$$
Lie_E g^{\alpha\beta} =(d-1) g^{\alpha\beta}, ~~ Lie_E
c_\gamma^{\alpha\beta} =(d-1) c_\gamma^{\alpha\beta}.
$$
The second statement is clear since the isomorphism
$$
(~,~): T^*M\to TM
$$
maps the unity $dt_1 = \eta_{1\alpha} dt^\alpha$ of the multiplication
(\ref{prod})
to the Euler vector field $E$. Proposition is proved.

\begin{cor}\label{cor1-2} Let $p^1(t)$, \dots, $p^n(t)$ be a system of flat
coordinates of the intersection form defined locally on
$M\setminus\Sigma$. Then there exists a function $F_*(p)$ such that
\beq\label{fstar}
{\partial^3 F_*(p)\over \partial p^i \partial p^j \partial p^k}
= G_{ia} G_{jb} {\partial t^\gamma\over \partial p^k} {\partial p^a\over
\partial t^\alpha} {\partial p^b \over \partial t^\beta}
c^{\alpha\beta}_\gamma(t).
\eeq
Here
$$
\left( G_{ij}\right) =\left( G^{ij}\right)^{-1}
$$
and a constant symmetric nondegenerate matrix $G^{ij}$ is defined by
$$
G^{ij} ={\partial p^i\over \partial t^\alpha} {\partial p^j\over \partial
t^\beta} g^{\alpha\beta}(t).
$$
The function $F_*(p)$ satisfies the following associativity equations
\beq\label{w-star}
{\partial^3 F_*(p)\over \partial p^i \partial p^j \partial p^a}
G^{ab} {\partial^3 F_*(p)\over \partial p^b \partial p^k \partial p^l}
=
{\partial^3 F_*(p)\over \partial p^l \partial p^j \partial p^a}
G^{ab} {\partial^3 F_*(p)\over \partial p^b \partial p^k \partial p^i},
~~i,\,j,\,k,\,l=1, \dots, n.
\eeq
For $d\neq 1$ $F_*(p)$ satisfies the homogeneity condition
\beq\label{q-h-star}
\sum_a p^a {\partial F_*\over \partial p^a} = 2\,F_* +{1\over 1-d} \sum
G_{ab} p^a p^b.
\eeq
\end{cor}

Observe that the definition of the function $F_*$ can be rewritten in the
following way
\beq\label{fstar1}
d\, \left( {\pal F_* \over \pal p^a \pal p^b}\right) = dp_a \cdot dp_b.
\eeq

We have constructed an {\it almost Frobenius structure} on the complement
$M_*=M\setminus\Sigma$. By definition, this means that all the axioms of Frobenius
manifold hold true but the constancy of the unity. We will called this object
{\it dual (almost) Frobenius manifold}. At the end of this Section we describe
the deformed flat coordinates of the dual Frobenius manifold. Denote $\nu$ the
parameter of the deformation (it was $z$ for the original Frobenius manifold).
By definition the deformed flat coordinates are independent 
flat functions $\tilde p(p;\nu)$ for the deformed flat connection  $\tilde\nabla_*$ defined on
$M_*\times {\mathbb C}$ by the formulae similar to (\ref{def-con1}). In the flat
coordinates $p^1$, \dots, $p^n$ for the intersection form they satisfy the
system
\beq\label{flat-star}
{\pal \xi_a\over \pal p^b} = \nu \cstar_{ab}^c(p) \xi_c, ~~ \xi_a ={\pal \tilde
p(p;\nu)\over \pal p^a}
\eeq
where
\beq\label{c-star}
\cstar_{ab}^c(p) := G^{cd} {\pal^3 F_*(p)\over \pal p^d \pal p^a \pal p^b}.
\eeq

{\bf Definition 8.} Any function $\tilde p=\tilde p(p;\nu)$ satisfying
(\ref{flat-star})
is called {\it twisted period} of the Frobenius manifold.

As usual one has an invariant pairing on twisted periods defined by the formula
\beq\label{pairing}
(\xi(\nu),\xi(-\nu)):= \xi_i(p;\nu) G^{ij}\xi_j(p;-\nu).
\eeq
The above expression does not depend on $p\in M_*$

Let us rewrite the system (\ref{flat-star}) in the original coordinates
$t^1$, \dots, $t^n$. 

\begin{prop} The gradients $\xi=(\xi_1, \dots, \xi_n)$ of 
the deformed flat coordinates on $M_*\times {\mathbb C}$
\beq\label{xi-dual}
\xi_\alpha:=\partial_\alpha \tilde p(t; \nu), ~~\alpha=1, \dots, n
\eeq
satisfy the following system of linear differential equations
\beq\label{g-m-dual}
\partial_\alpha\xi\cdot {\mathcal U} = \xi\cdot \left( {\mathcal V}+\nu-{1\over 2}\right) \,
C_\alpha.
\eeq
\end{prop}

The notations are the same as above.

\pf By definition the differential equations for the deformed flat
coordinates have the form
$$
\partial_\alpha\xi_\beta -\Gamma_{\alpha\beta}^\gamma\xi_\gamma -\nu\, \hat
c_{\alpha\beta}^\gamma\xi_\gamma=0.
$$
Here $\Gamma_{\alpha\beta}^\gamma$ are the Christoffel coefficients of the
Levi-Civita connection for the metric \ref{int-form}), $\hat c_{\alpha\beta}^\gamma$
are the structure constants of the dual Frobenius manifold. Multiplying
this equation by ${\mathcal U}_\rho^\alpha$ and using (\ref{christ}) 
and the definition (\ref{star})
of the multiplication law on the dual Frobenius manifold we arrive at
(\ref{g-m-dual}). \epf

For $\nu=0$ we obtain the original periods
$$
\tilde p(p;0)\equiv p
$$
i.e., the flat functions for the metric $(~,~)$.
The twisted periods 
\beq\label{odd}
\varpi(t)=\tilde p(t;\nu=1/2)
\eeq
will be called {\it odd periods}. This name is
motivated by the following construction.

Let us introduce a Poisson bracket on the Frobenius manifold by the formula
\beq\label{pb}
\{ f, g\} := <df, {\mathcal V}\,dg>.
\eeq
The flat coordinates $t^1$, \dots, $t^n$ are Darboux coordinates for this
Poisson bracket:
\beq\label{darb}
\{ t^\alpha, t^\beta\} = \eta^{\alpha\gamma} {\mathcal V}^\beta_\gamma.
\eeq
The Poisson structure does not degenerate {\it iff} the operator ${\mathcal V}$ is
invertible.

\begin{prop}\label{prop-pb} Let $\varpi_1(t)$, $\varpi_2(t)$  be any 
two odd periods. Then their Poisson bracket (\ref{pb}) is constant.
\end{prop}

Proof is given by a straightforward calculation of the derivatives
$\partial_\alpha\{ \varpi_1, \varpi_2\}$ using (\ref{g-m-dual}) for $\nu={1\over
2}$.

\begin{cor} A system of $n$ independent odd periods $\varpi_1(t)$, 
\dots, $\varpi_n(t)$
 gives another system of Darboux coordinates for the
Poisson bracket (\ref{pb}).
\end{cor}

The generating function of the canonical transformation 
$$
\left( t^1, \dots, t^n\right) \mapsto \left( \varpi_1(t), \dots,
\varpi_n(t)\right)
$$
between the two systems
of Darboux coordinates will play an important role in the constructions of
Section 4.

Let us now consider the dependence of the twisted periods on $\nu$. It turns out
that, instead of the differential equation with rational coefficients given by
(\ref{def-con2}) the crucial role in the theory of almost Frobenius manifolds
plays a {\it difference equation in} $\nu$.

\begin{lemma} If $\xi=\xi(t;\nu)$ is a solution to (\ref{g-m-dual}) then $\partial_1\xi$
is a solution to the same system with a shift $\nu\mapsto \nu-1$.
\end{lemma}

We will represent the claim of the lemma in the form
\beq\label{shift}
\partial_1\xi(t;\nu) =\xi(t;\nu-1)
\eeq
and, respectively for the twisted periods
\beq\label{p-shift}
{\pal\over\pal t^1} \tilde p(t;\nu) = \tilde p(t;\nu-1).
\eeq

\pf Dependence of $\xi(t;\nu)$ on $t^1$ is determined, due to (\ref{g-m-dual}), by the
system
\beq\label{ur1}
\dot \xi\cdot {\mathcal U} = \xi\cdot \left( {\mathcal V}+\nu-{1\over 2}\right), ~~\dot \xi:=
\partial_1\xi. 
\eeq
So one can rewrite the equation (\ref{g-m-dual}) as
$$
\pal_\alpha\xi\cdot {\mathcal U} = \dot\xi \cdot C_\alpha {\mathcal U}
$$
where the commutativity ${\mathcal U} C_\alpha = C_\alpha {\mathcal U}$ of operators of
multiplication has been used. Because of invertibility of the operator ${\bf
U}$ on $M_*$ one has
$$
\pal_\alpha \xi = \dot \xi \cdot C_\alpha.
$$
Differentiation (\ref{g-m-dual}) w.r.t. $t^1$ and using $\partial_1 {\mathcal U}(t)={\rm id}$,
$\partial_1 C_\alpha =0$ yields
$$
\pal_\al \dot \xi \cdot {\mathcal U} + \pal_\al \xi = \dot\xi \cdot \left( {\mathcal V} + \nu
-{1\over 2}\right) C_\al.
$$
Using the above equation we arrive at
$$
\pal_\al \dot\xi =\dot\xi \cdot \left( {\mathcal V} + \nu
-{3\over 2}\right) C_\al.
$$
This equation coincides with (\ref{g-m-dual}) up to shift $\nu\mapsto \nu-1$.
The Lemma is proved.\epf

Using (\ref{g-m-dual}) one can rewrite the shift equation (\ref{shift}) in the
form
\beq\label{ab-shift}
\xi(t;\nu-1) = \xi(t;\nu) \left[ {\mathcal A}(t) \nu + {\mathcal B}(t)\right], ~~
{\mathcal A}(t) = {\mathcal U}^{-1}, ~~ {\mathcal B}(t) = \left( {\mathcal V} + \nu -{1\over 2}
\right){\mathcal U}^{-1}.
\eeq
Observe that the shift operator (\ref{ab-shift}) {\it changes the sign} of the
invariant bilinear form (\ref{pairing}): if $\xi(t;\nu -1)$ and $\xi(t;\nu)$
are related by (\ref{ab-shift}) then
\beq\label{sign}
(\xi(\nu-1), \xi(-\nu+1)) = - (\xi(\nu), \xi(-\nu)).
\eeq

We are now ready to give a precise definition of almost Frobenius manifold
and prove a Reconstruction Theorem inverting the above ``almost duality''
$$
M\mapsto M_*.
$$ 
For simplicity we will only consider the
case $d\neq 1$. Moreover we will present the definition of almost Frobenius
manifolds only in the flat coordinates leaving the coordinate-free formulation,
along the lines of Section \ref{crash}, as an exercise for the reader.

{\bf Definition 9.} {\it Almost Frobenius structure} of the charge $d\neq 1$
on the manifold $M_*$ is a structure of a Frobenius algebra on the tangent
planes $T_pM_* = \left(~*~,(~,~)_p\right)$ depending (smoothly, analytically etc.)
on the point $p\in M_*$. It must satisfy the following axioms.

{\bf AFM1.} The metric $(~,~)_p$ is flat.

{\bf AFM2.} In the flat coordinates $p^1$, \dots, $p^n$ for the metric,
\beq\label{gram}
(dp^i, dp^j) = G^{ij}
\eeq
the structure constants of the multiplication
\beq\label{mult}
{\pal\over\pal p^i} * {\pal\over \pal p^j} = \cstar_{ij}^k(p) {\pal\over\pal
p^k}
\eeq
can be locally represented in the form
\beq\label{3der}
\cstar_{ij}^k(p)=G^{kl}{\pal^3 F_*(p)\over \pal p^l\pal p^i\pal p^j}
\eeq
for some function $F_*(p)$. The function satisfies the homogeneity equation
\beq\label{homo}
\sum_{i=1}^n p^i {\pal F_*(p)\over \pal p^i} = 2 \, F_*(p) + {1\over 1-d} (p,p).
\eeq
The Euler vector field
\beq\label{E}
E={1-d\over 2} \sum_{i=1}^n p^i {\pal \over \pal p^i}
\eeq
is the unity of the Frobenius algebra.

{\bf AFM3.} There exists a vector field 
\beq\label{e}
e=e^k(p) {\pal\over\pal p^k}
\eeq
being an invertible element of the Frobenius algebra $T_pM_*$ for every $p\in
M_*$ such that the operator
$$
\tilde p \mapsto e\tilde p
$$
acts by a shift $\nu\mapsto \nu -1$ on the solutions of the deformed flat
connection equations
\beq\label{def-star}
{\pal^2 \tilde p\over \pal p^i \pal p^j} = \nu \, \cstar_{ij}^k(p) {\pal \tilde
p\over \pal p^k}.
\eeq

\begin{theorem}\label{reconstruct} Let us define on an almost Frobenius manifold.
$M_*$ a new multiplication $\cdot$ and a new bilinear form $<~,~>$ by 
\beq\label{dot}
u\cdot v := {u*v\over e}
\eeq
\beq\label{angle}
<u,v> := \left( u, {v\over e}\right).
\eeq
These multiplication and bilinear form gives a Frobenius structure of the charge
$d$ on $M_*$
with the unity $e$ and the Euler vector field $E$.
\end{theorem}

\pf Using (\ref{def-star}) we can rewrite the shift operator $\nu\mapsto\nu-1$
as follows
\beq\label{pab-shift}
\xi(p; \nu-1) =\xi(p; \nu)(A(p)\nu + B(p)), ~~\xi(p;\nu) = \left(
\xi_i(p;\nu)\right), ~~\xi_i(p;\nu)={\pal \tilde p(p;\nu)\over \pal p^i}
\eeq
(cf. (\ref{ab-shift}) above) where
\eqa\label{ab}
&&
A(p) = (A^i_j(p)), ~~A^i_j(p) = e^k(p) \cstar_{kj}^i(p)
\nn\\
&&
B(p)=(B^i_j(p)), ~~ B^i_j(p) = \pal_j e^i(p).
\eeqa
Since the operator (\ref{pab-shift}) changes the sign of the pairing
(\ref{pairing}), the operators $A=A(p)$, $B=B(p)$ satisfy
$$
(Au, v)=(u, Av), ~~(Bu,v)+(u,Bv) =-(u,Av)
$$
for any two tangent vectors $u$, $v$. The second of these two equations implies
\beq\label{de}
\pal^i e^j + \pal^j e^i = - e^k \cstar_k^{ij}
\eeq
(all raising and lowering of Latin indices in this Section is done by means
of the matrix $\left(G^{ij}\right)$ and its inverse $\left(G_{ij}\right)$).

Let us first prove that the metric $<~,~>$ is flat. It is convenient to work
with the contravariant metric on $T^*M_*$. The Gram matrix of it in the
coordinates $p^i$ is given by
\beq\label{eta}
\eta^{ij}(p):=<dp^i, dp^j> = e^k(p) \cstar^{ij}_k(p).
\eeq
Let us compute the Christoffel coefficients of the Levi-Civita connection for
the metric $<~,~>$. Let us denote these coefficients $\gamma_{ij}^k(p)$; we put
\beq\label{ctvr}
\gamma^{ij}_k(p):= -\eta^{is}(p) \gamma_{sk}^j(p)
\eeq
(cf. (\ref{christ0})).

\begin{lemma}\label{lemma1} The Christoffel coefficients (\ref{ctvr}) of the
Levi-Civita connection for the metric (\ref{eta}) are given by one of the
following two equivalent formulae
\eqa\label{gamma1}
&&
\gamma^{ij}_k=-\pal^i\pal_k e^j
\\
&&
\gamma^{ij}_k= \cstar^i_{ks}\pal^s e^j.\label{gamma2}
\eeqa
\end{lemma}

\pf Let us first check equality of (\ref{gamma1}) and (\ref{gamma2}). Indeed,
from
(\ref{p-shift}) it follows that the functions $\tilde p=e^1(p)$, \dots, $\tilde
p=e^n(p)$ give a
basis in the space of solutions of the equation (\ref{def-star}) for $\nu=-1$.
Spelling this equation out yields the needed equality. 

We are now to check that the Christoffel coefficients $\gamma^{ij}_k$ satisfy
the following two equations (see \cite{cime}, equations (3.26) and (3.27))
uniquely determining the Levi-Civita connection
\eqa\label{cime1}
&&
\pal_k\eta^{ij} = \gamma^{ij}_k+\gamma^{ji}_k   
\\
&&
\eta^{is}\gamma{jk}_s = \eta^{js}\gamma^{ik}_s.
\label{cime2}
\eeqa
Indeed, the first of these two equations follows from (\ref{gamma1}) and from
the identity
\beq\label{eta1}
\eta^{ij}=-\pal^ie^j -\pal^je^i
\eeq
equivalent to 
(\ref{de}). The second one follows from (\ref{gamma2}) and from associativity
of the multiplication $*$. The Lemma is proved. \epf

We will now prove that the metric (\ref{eta}) is flat. Because of (\ref{cime2})
the connection is torsion-free. It remains to prove

\begin{lemma}\label{lemma2} The curvature of the connection (\ref{gamma1}),
(\ref{gamma2}) vanishes.
\end{lemma}

\pf We are to prove the following identity (see \cite{cime}, formula (3.33))
$$
\eta^{is} \left( \pal_s \gamma_l^{jk} - \pal_l \gamma_s^{jk}\right) +
\gamma^{ij}_s \gamma^{sk}_l - \gamma^{ik}_s \gamma^{si}_l=0.
$$
To this end we are to substitute the formula (\ref{gamma1}) in the first
bracket and also the expressions (\ref{eta}) and (\ref{gamma2}) in the remaining
terms. After such substitution the first bracket will be vanishing because of 
equality of mixed derivatives, and the remaining terms will be equal to zero
because of associativity of the product $*$. The Lemma is proved. \epf

We will now proceed to studying the properties of the multiplication
(\ref{dot}).

\begin{lemma}\label{lemma5} The multiplication (\ref{dot}) is commutative and
associative. The vector $e$ is the unity of it. The bilinear form $<~,~>$
is symmetric nondegenerate and invariant w.r.t. the multiplication.
\end{lemma}

\pf Rewriting the formulae (\ref{dot}) and (\ref{angle}) in the form
$$
u\cdot v = u*v*e^{-1}
$$
$$
<u,v> = (u*v, e^{-1})
$$
we obtain
$$
<u\cdot v, w> = (u*v*w, e^{-2}).
$$
The statements of the Lemma easily follow from the above formulae. \epf

We are now to check the main property {\bf FM2} of the multiplication
(\ref{dot}). Let us denote $\nabla_k$ the covariant derivatives w.r.t. the
Levi-Civita connection for the metric $<~,~>$. Introduce the operators
$$
\nabla^k := \eta^{ks}\nabla_s.
$$
It suffices to prove 

\begin{lemma}\label{lemma6} The coefficients $c^{ij}_k = c^{ij}_k(p)$ defined by
$$
dp^i\cdot dp^j = c^{ij}_k(p) dp^k
$$
satisfy
\beq\label{nabla-c}
\nabla^l c^{ij}_k = \nabla^i c^{lj}_k.
\eeq
\end{lemma}

\pf The main trick in the proof is the coincidence of the two multiplications
$*$ and $\cdot$ on the {\it cotangent} spaces $T^*M_*$. So $c^{ij}_k$ in 
(\ref{nabla-c}) can be replaced by $\cstar^{ij}_k$. So we are to prove symmetry
in $i$, $j$, $l$ of the covariant derivative
$$
\nabla^l \cstar^{ij}_k =\eta^{ls} \pal_s \cstar^{ij}_k -\gamma^{li}_q
\cstar^{qj}_k - \gamma^{lj}_q \cstar ^{iq}_k + \gamma^{lq}_k \cstar^{ij}_q.
$$
The derivative $\pal_s \cstar^{ij}_k$  is symmetric in $s$, $k$ due to
(\ref{3der}). Using this symmetry and also the equation (\ref{cime1}) we recast
the r.h.s. of the last expression into the form
$$
= \pal_k \left( \eta^{ls}\cstar_s^{ij}\right) - \gamma^{li}_q \cstar^{qj}_k
-\gamma^{lj}_q \cstar^{qi}_k - \gamma^{ql}_k \cstar^{ij}_q.
$$
In the first term we can replace back $\cstar_s^{ij}$ by $c_s^{ij}$. After such a
replacement the symmetry of the first term in $i$, $j$, $l$ becomes evident.
The remaining three terms can be rewritten, using (\ref{gamma2}) as the $k$-th
component of the following one-form
$$
- \left( dp^l * dp^j * de^i + dp^l * dp^i * de^j + dp^i * dp^j * de^l\right).
$$
The needed symmetry becomes obvious. The Lemma is proved.\epf

\begin{lemma}\label{lemma3} The unity vector field $e$ satisfies (\ref{e-flat}).
\end{lemma}

\pf Using the formula (\ref{eta}) and also the expression (\ref{gamma2}) for the
Christoffel coefficients we compute
$$
\nabla^i e^k = \eta^{is} \pal_s e^k - \gamma^{ik}_s e^s
=e^l \cstar^{is}_l \pal_s e^k - \cstar^i_{sl} \pal^l e^k \, e^s =0.
$$
The Lemma is proved. \epf

It remains to settle the quasihomogeneity property {\bf FM3}. To this end the
following quasihomogeneity of the deformed flat coordinates $\tilde p(p; \nu)$
will be useful
\beq\label{quasi-p}
Lie_E \tilde p(p; \nu)= \left[{1-d\over 2} + \nu\right] \, \tilde p(p; \nu).
\eeq 
In particular, the components $e^k(p)$ of the vector field $e$ satisfy
\beq\label{quasi-e}
Lie_E e^k = -{1+d\over 2} e^k, ~~k=1, \dots, n
\eeq
(see the proof of Lemma \ref{lemma1}). 

\begin{lemma}\label{lemma4} The vector field $E$ is linear in the flat
coordinates for the metric $<~,~>$. It satisfies 
\beq\label{comm}
[e,E]=e.
\eeq
\end{lemma}

\pf Using (\ref{gamma2}) and also 
\beq\label{e-unity}
E^k \cstar_{kj}^i = \delta_j^i
\eeq
(see Axiom {\bf AFM2}) we obtain the following expression for the first
covariant derivative of the vector field $E$
$$
\nabla^i E^j = {1-d\over 2} \eta^{ij} - \pal^i e^j.
$$
For the second covariant derivative $\nabla^k\nabla^i E^j$ we obtain, using
(\ref{cime1}) and also equality of (\ref{gamma1}) and (\ref{gamma2})
$$
\nabla^k\nabla^i E^j= {1-d\over 2} (\gamma^{ij}_s+\gamma^{ji}_s) + \eta^{ks}
\gamma_s^{ij} -\gamma_q^{ki}\left[ {1-d\over 2} \eta^{qj}- \pal^q e^j\right]
-\gamma^{kj}_q \left[ {1-d\over 2} \eta^{iq} - \pal^i e^q\right].
$$
After substitution of the expression (\ref{gamma2}) for the Christoffels we
rewrite the last equation as follows
$$
= {3-d\over 2} \eta^{ks} \cstar^{iq}_s \pal_q e^j +{1-d\over 2} \eta^{ks}
\cstar^{jq}_s \pal_q e^i
-{1-d\over 2} \eta^{js} \cstar^{kq}_s \pal_q e^i - {1-d\over 2} \eta^{is}
\cstar^{kq}_s \pal_q e^j
$$
$$
+ \cstar^k_{qs} \pal^s e^i \pal^q e^j + \cstar^k_{qs} \pal^s e^j \pal^i e^q.
$$
Replacing $\cstar^{iq}_s$, $\cstar^{jq}_s$, $\cstar^{kq}_s$, and $\cstar^{kq}_s$
by $c^{iq}_s$, $c^{jq}_s$, $c^{kq}_s$, and $c^{kq}_s$ resp. and using symmetries
like
$$
\eta^{ks} c^{jq}_s = \eta^{js} c^{kq}_s
$$
etc. we arrive at the following expression
$$
= \eta^{is} c^{kq}_s \pal_q e^j + c^{qk}_s \left( \pal^i e^s + \pal^s
e^i\right) \pal_q e^j.
$$
The expression in the parenthesis in the second term is equal to $-\eta^{is}$, 
due to (\ref{de}). Therefore $\nabla^k\nabla^i E^j=0$.

To prove (\ref{comm}) we are to use (\ref{quasi-e}). So,
$$
[e,E]^k = e^i \pal_i {1-d\over 2} p^k - E^i\pal_i e^k = {1-d\over 2} e^k +
{1+d\over 2} e^k = e^k.
$$
The Lemma is proved. \epf

The last step in the proof of the Theorem \ref{reconstruct} is

\begin{lemma} The linear vector field $E$ satisfies equations of Axiom {\bf
FM3}.
\end{lemma}

\pf Let us first prove that
$$
Lie_E \eta^{ij} \equiv E^k\pal_k \eta^{ij} -\pal_k E^i \eta^{kj} - \pal_k E^j
\eta^{ik} = (d-2) \eta^{ij}.
$$
Using (\ref{cime1}) we rewrite the last equation in the form
$$
= E^k (\gamma^{ij}_k + \gamma^{ji}_k) - (1-d) \eta^{ij}.
$$
Substituting the expression (\ref{gamma2}) for the Christoffels and using
(\ref{e-unity}) and (\ref{de}) we arrive at
$$
= \pal^i e^j + \pal^j e^i - (1-d) \eta^{ij} = (d-2) \eta^{ij}.
$$

Next, we are to prove that
$$
Lie_E c^{ij}_k  \equiv E^s\pal_s c^{ij}_k -\pal_s E^i c^{sj}_k -\pal_s E^j
c^{is}_k + \pal_k E^s c^{ij}_k = (d-1) c^{ij}_k.
$$
The last three terms give $-(1-d) c^{ij}_k$. In the first term, because of 
equality $c^{ij}_k= \cstar^{ij}_k$ and because of (\ref{3der}) we can
interchange the indices $s$ and $k$. After this we rewrite this term as follows
$$
E^s \pal_k c^{ij}_s = \pal^k \left( E^s \cstar^{ij}_s\right) - \pal_k E^s
c^{ij}_s = - c^{ij}_k.
$$
Here we use again (\ref{e-unity}). This proves the Lemma and also the Theorem
\ref{reconstruct}. \epf

{\bf Example 1.} Let $\left( A={\rm span}(e_1, \dots, e_n), <~,~>\right)$
be a Frobenius algebra with the trivial grading
$$
q_\alpha={\rm deg}\, e_\alpha =0, ~~\alpha=1, \dots, n.
$$
It carries a structure of a trivial Frobenius manifold $M=A$ with
$$
F(t)={1\over 6}< t^3, e>, ~~t=t^\alpha e_\alpha\in A.
$$
Here $e\in A$ is the unity. The Euler vector field is equal to $E=t^\alpha
\partial_\alpha$. The dual (almost) Frobenius manifold $M^*$ can be identified
with the set of all invertible elements $x=x^\alpha e_\alpha$ of $A$ with the
metric
$$
\left( {\partial\over \partial x^\alpha}, {\partial\over \partial
x^\beta}\right) =  <e_\alpha,e_\beta>
$$
and with the multiplication defined by
$$
\left( {\partial\over \partial x^\alpha}* {\partial\over \partial
x^\beta},  {\partial\over \partial x^\gamma}\right) = 2 (e_\alpha \cdot
e_\beta\cdot e_\gamma, x^{-1}).
$$
Integrating one obtains
\beq\label{star-triv}
F_*(x) ={1\over 2} \left( x^2, \log x^2\right).
\eeq
The map
\eqa\label{miura}
&&
M^* \to M
\nn\\
&&
x \mapsto t={1\over 4} x^2
\eeqa
transforms the metric $(~,~)$ to the intersection form of the trivial Frobenius
manifold $M$. The latter can be recast into the following bilinear form on
$T^*M$ depending linearly on the coordinates
$$
\left( dt^\alpha, dt^\beta\right) = c^{\alpha\beta}_\gamma t^\gamma.
$$
This metric and the flat coordinates (\ref{miura}) for it was first considered by
A.Balinski and S.P.Novikov \cite{bn} in their theory of linear Poisson brackets of
hydrodynamic type. The solution (\ref{star-triv}) to the equations of associativity was found
in Appendix to \cite{pisa, cmp92} (in a different but equivalent form). 

We will end this Section with a slightly more general construction of the almost
dual Frobenius manifold. Let us introduce a one-parameter family of Frobenius
algebra structures on the tangent spaces $T_tM$ by
\eqa\label{lam-dual}
&&
u *_\lambda v := {u\cdot v\over E-\lambda e}
\nn\\
&&
(u,v)_\lambda := <u, {v\over E-\lambda e}>.
\eeqa
It is easy to see that, for any $\lambda$ the above formulae define on $M_*(\lambda):= M\setminus
\Sigma_\lambda$ a structure of almost Frobenius manifold. For $\lambda=0$
one obtains the old definition $M_*(0)=M_*$. For $\lambda\to\infty$ after a
suitable rescaling the
Frobenius structure (\ref{lam-dual}) goes to the original Frobenius structure.

\setcounter{equation}{0}
\setcounter{theorem}{0}

\section{Twisted period mapping and its monodromy}\label{sec3}

In this section we will describe the analytic properties of the twisted period
mapping for an arbitrary semisimple Frobenius manifold.

Let $\nu$ be an arbitrary complex number,
$$
q=e^{2\pi i\nu}.
$$
We first introduce one more argument $\lambda$ of the twisted periods
doing a shift 
$$
\tilde p(t;\nu) \mapsto e^{-\lambda \partial_1} \tilde p(t;\nu) =
\tilde p(t^1-\lambda, t^2, \dots, t^n; \nu)
$$
The gradients $\xi_\alpha=\partial_\alpha\tilde p$ of these functions can be found from the following system of
linear differential equations.

\begin{lemma}\label{lemma2-1} Near any point $t\in M\setminus\Sigma_\lambda$ there
exist $n$ independent functions $\tilde p^1(\nu;\lambda;t)$, \dots,
$\tilde p^n(\nu;\lambda;t)$ such that their gradients $\xi=(\xi_1, \dots, \xi_n)$,
$\xi_\alpha=\partial_\alpha \tilde p$, satisfy the system
\eqa\label{g-m-l-mu}
&&
\partial_\alpha\xi \cdot ({\mathcal U}-\lambda) = \xi\,\left({\mathcal V}+\nu -{1\over
2}\right)\, C_\alpha
\\
&&
\partial_\lambda \,\xi \cdot  ({\mathcal U}-\lambda) = \xi\,\left(-{\mathcal V}-\nu
+{1\over      
2}\right).\label{g-m-l-lam}
\eeqa
\end{lemma}

To prove Lemma it suffices to check compatibility of the system.
We leave it as an exercise to the reader.

Observe that for $\nu=0$ the system (\ref{g-m-l-mu}), (\ref{g-m-l-lam}) coincides 
with (\ref{g-m0}), (\ref{g-m01}). More generally, it coincides with the
equations defining deformed flat coordinates on the almost dual Frobenius
manifold $M_*(\lambda)$ (see (\ref{lam-dual}) above).

Our nearest goal is to describe the monodromy of solutions of the 
system (\ref{g-m-l-mu}), (\ref{g-m-l-lam}).
This will be done for an arbitrary semisimple Frobenius manifold.

We first rewrite the system (\ref{g-m-l-mu}), (\ref{g-m-l-lam}) for the twisted
periods in the canonical coordinates. We will denote $\tilde p(\nu;
\lambda; u)$ the function $\tilde p(\nu;\lambda;t(u))$ when it cannot lead to a
confusion.

Let $\phi_i(\nu;\lambda;t)$ be the components of the one-form $\xi=d\tilde p$
in the moving frame $f_1$, \dots, $f_n$,
\beq\label{p-phi}
\phi_i ={\partial_i\tilde p\over \psi_{i1}}.
\eeq

\begin{lemma}\label{lemma2-2} In the new coordinates the system 
(\ref{g-m-l-mu}), (\ref{g-m-l-lam}) reads
\eqa\label{system1}
&&
\partial_i\phi = \left( {E_i \left({1\over 2}-\nu +V\right)\over
\lambda-u_i}+V_i\right)\,\phi
\\
&&
(U-\lambda) {d\phi\over d\lambda} = \left({1\over 2}-\nu +V\right)\,\phi
\label{system2}
\eeqa
(we write now $\phi$ as a column-vector). Here $E_i$ is the matrix having
all the entries zero but $(E_i)_{ii}=1$, the skew-symmetric matrix $V_i$
has the form
$$
V_i ={\rm ad}_{E_i} {\rm ad}_U^{-1} (V).
$$
\end{lemma}

Proof is similar to Proposition H.2 in \cite{cime}. We omit it.

If $\phi$ is a solution to the system (\ref{system1}), (\ref{system2}) then $d\phi/d\lambda$ is 
a solution of the same system with $\nu\mapsto \nu-1$. So it suffices to
describe the monodromy of the system assuming that
$$
{\rm Re}\, \nu > -{1\over 2}.
$$

Let us cover the Frobenius manifold $M_{s\, s}$ with convenient charts.
We choose a real number $0\leq \varphi <2\pi$ and we define an open subset
$M_{s\, s}^0 \subset M_{s\, s}$ containing all the points $t\in M_{s\, s}$
such that their canonical coordinates $u_1(t)$, \dots, $u_n(t)$ satisfy the
following condition: the rays $L_1$, \dots, $L_n$ on the complex plane of the
form
\beq\label{cuts}
L_j =\left\{ u_j + i\, \rho\, e^{-i\varphi} ~|~ 0\leq \rho <\infty\right\} ,
~~j=1, \dots, n
\eeq
must not intersect. On $M_{s\, s}^0$ we can order the canonical
coordinates $u_1$, \dots, $u_n$ in such a way that the rays $L_1$, \dots,
$L_n$ exit from the infinite point of the complex plane in the 
counter-clockwise
order. After this ordering we are able to define the matrix-valued
functions $\Psi(u)$ and $V(u)$ as it was explained above. We can define
therefore the linear differential operator with rational coefficients
\beq\label{oper}
{d\over dz} -\left( U+{V\over z}\right)
\eeq
and compute it monodromy data $\left( \hat\mu, e,  R, S, C\right)$ at each point of 
$M_{s\, s}^0$. Here the $n\times n$ matrices $\hat\mu$, $R$ describe monodromy
at the origin, $e$ is an eigenvector of the matrix $\hat\mu$, $S$ is the Stokes
matrix of the operator (\ref{oper}) computed with respect to the line 
\beq\label{ell}
\ell =\{ \arg \, z=\varphi\}
\eeq
with its natural orientation, $C$ is the central connection matrix (see the
definitions and the full list of constraints for the monodromy data in
\cite{pain}). 
Observe that $S$ is an upper triangular matrix due to the
above choice of ordering of the entries $u_1$, \dots, $u_n$ of the diagonal
matrix $U$. The central result of the theory of semisimple Frobenius manifolds
says that the monodromy data are constant on every connected piece of $M_{s\,
s}^0$. The Frobenius manifold structure on any such a piece 
can be reconstructed by an algebraic
procedure starting from the solution of a suitable Riemann - Hilbert
boundary value problem  with the boundary conditions given in terms of the
monodromy data. We will denote $Fr(\hat\mu, R, e, S, C)$
such a Frobenius structure on any
connected component of $M_{s\, s}^0$  characterized by the monodromy data 
$(\hat\mu, R, e, S, C)$. 
As a consequence of the general theory of Riemann - Hilbert problems we derive 
that the image of the map
$$
M_{s\, s}^0 \to {\mathbb C}^n \setminus \diag , ~~t \mapsto (u_1(t), \dots, u_n(t))
$$
is a complement to a closed analytic subset in ${\mathbb C}^n$. The gluing of the 
patches $Fr(\hat\mu, R, e, S,$ $ C)$ 
along the boundaries
$$
\arg \, (u_i -u_j) ={\pi\over 2} -\varphi ~~{\rm for ~some}~ i\neq j
$$
is given by an action of the braid group
$B_n$ on the monodromy data described in \cite{cime, pain}.

In every patch $Fr(\hat\mu, R, e, S, C)$ we will construct a fundamental matrix of
solutions
of the Fuchsian system (\ref{fuchs1}) depending on $u_1$, \dots, $u_n$ according to the 
equations (\ref{fuchs2}). 

\begin{theorem}\label{lef} Let $q:= e^{2\pi i \nu}$ be not a root of the characteristic
equation
\beq\label{ch-q}
\det \left( q\, S + S^T\right) =0.
\eeq
Then there exist $n$ linearly independent solutions $\phi^{(1)}$, 
\dots, $\phi^{(n)}$ of the system (\ref{fuchs1}) analytic in $(\lambda, u_1, \dots,
u_n)$ on 
\beq\label{dom}
{\mathbb C} \setminus \cup_j L_j \, \times \, Fr(\hat\mu, R, e, S, C)
\eeq
such that the monodromy transformations $M_1$, \dots, $M_n$ along the small
loops encircling counter-clockwise the points $u_1$, \dots, $u_n$ are
reflections
\beq\label{reflect}
M_i \phi^{(j)} = \phi^{(j)} - q^{1/2}\left( \check\phi^{(i)},
\phi^{(j)}\right)_q \phi^{(i)}
\eeq
w.r.t.
the bilinear form 
\beq\label{q-form}
\left( \check\phi^{(i)}, \phi^{(j)}\right)_q := \left( q^{1/2} S+q^{-1/2}S^T\right)_{ij}.
\eeq
\end{theorem}

\pf For ${\rm Re}\, \nu << 0$ we can construct a fundamental matrix
$\Phi(\lambda)$ of solutions to (\ref{system2}) applying a Laplace-type transform to
a fundamental matrix $Y(z)$ of solutions to (\ref{isomono}):
\beq\label{lap0}
\Phi(\lambda) ={i\over \sqrt{2\pi}} (1+q^{-1}) 
\int_0^\infty Y(z) e^{-\lambda\,z} {dz\over z^{\nu+{1\over 2}}}.
\eeq
Technically it is more convenient, following \cite{balser,cime, pain}, to use a sort of inverse
transform expressing solutions to (\ref{isomono}) via Laplace-type integrals applied to the
solutions to the system (\ref{system2}).

We rewrite the Fuchsian system (\ref{system2}) in the standard way
\beq\label{fuchs1}
{d\phi\over d\lambda} = \sum_{i=1}^n {A_i \over \lambda -u_i} \,\phi
\eeq
where
\beq\label{fuchs2}
A_i := E_i \, \left( \nu -{1\over 2} -V\right), ~~i=1, \dots, n.
\eeq
The matrices $E_i$ were defined in (\ref{entr}). This system has Fuchsian singularities 
at the points $\lambda=u_1$, \dots, $\lambda=u_n$, $\lambda=\infty$. The
dependence on $u_1$, \dots, $u_n$ imposed by (\ref{system1}) reads
\beq\label{fuchs3}
{\pal \phi\over \pal u_i} = \left( V_i - {A_i\over \lambda-u_i}\right)\, \phi.
\eeq

We define
the needed solutions $\phi^{(1)} , \dots, \phi^{(n)}$ according to their
behaviour near the finite singularities of the system.
 
Let us first consider the generic case $\nu$ being not a half-integer.

\begin{lemma} For 
\beq\label{half}
\nu\not\in{\mathbb Z}+{1\over2} 
\eeq
there exist $n$ solutions
of the Fuchsian system (\ref{fuchs1}) analytic in
\beq\label{dom-cuts}
\lambda \in {\mathbb C}\setminus \cup_j L_j
\eeq
of the form
\beq\label{local}
\phi^{(i)} ={\sqrt{2\pi}\over \Gamma\left({1\over 2}+\nu\right)}(u_i-\lambda)^{\nu -{1\over 2}} \left[ e_i +O(u_i-\lambda)\right],
~~\lambda\to u_i
\eeq
where the expression in the square brackets is analytic in a small neighborhood
of $\lambda=u_i$. 
Here $e_i$ is the column-vector having the i-th component 1 and all other
components zero.
The solutions are determined uniquelly after chosing
of a branch of the functions $(u_i-\lambda)^{\nu -{1\over 2}}$.
\end{lemma}

We choose the normalizing factor on (\ref{local}) in order to meet the shift
condition (\ref{shift}).

\pf The $n-1$ eigenvalues of the matrix $A_i$  are zeroes and one
eigenvalue is equal to $\nu-1/2$. So the eigenvalues of the monodromy
transformation $M_i$ are all equal to 1 but a simple eigenvalue $-q$. 
In the nonresonant case 
$$
\nu-{1\over 2} \not\in {\mathbb Z}
$$
the unique eigenvector of $M_i$ with the eigenvalue $-q$ can be represented,
after an appropriate normalization, in the form (\ref{local}). Lemma is proved.

Observe that any solution $\phi$ of the Fuchsian system (\ref{fuchs1}) near $\lambda=u_i$
can be uniquelly
represented in the form
$$
\phi = g\, \phi^{(i)} + analytic
$$
for some constant $g$. Here {\it analytic} is a solution of the same system
analytic at $\lambda=u_i$. Particularly, one can find a matrix of constants
$G=\left(g^{ij}\right)$ such that
\beq\label{gij}
\phi^{(i)} = g^{ji} \phi^{(j)} + analytic ~~{\rm near} ~\lambda=u_j.
\eeq
By definition we have
$$
g^{ii}=1, ~~i=1, \dots, n.
$$

\begin{lemma} \label{lemma35} The action of the monodromy transformation $M_j$ onto the
solution $\phi^{(i)}$ is given by the formula
\beq\label{monoij}
M_j \phi^{(i)} = \phi^{(i)} - (q+1) g^{ji} \phi^{(j)}.
\eeq
\end{lemma}

\pf From (\ref{local}) we have
$$
M_i\phi^{(i)} =-q\, \phi^{(i)}.
$$
This gives (\ref{monoij}) for $j=i$. If
$$
\phi^{(i)} =g^{ji} \phi^{(j)} +\psi
$$
with $\psi$ analytic near $\lambda=u_j$ then
$$
M_j \phi^{(i)}  = M_j \left( g^{ji} \phi^{(j)} + \psi \right)
$$
$$
=-q\, g^{ji} \phi^{(j)} + \psi = \phi^{(i)} - (q+1) g^{ji} \phi^{(j)}.
$$
Lemma is proved.

Similarly, continuing the solution $\phi^{(i)}$ clockwise around $\lambda=u_j$
we obtain the transformation
\beq\label{anti-mono}
M_j^{-1}\phi^{(i)}  = \phi^{(i)} -(1+q^{-1}) g^{ji} \phi^{(j)}.
\eeq

Let us introduce another system of branch cuts $\hat L_1$, \dots, $\hat L_n$
on the complex $\lambda$-plane opposite to $L_1$, \dots, $L_n$ resp. We
introduce solutions $\hat \phi^{(1)}$, \dots, $\hat \phi^{(n)}$ to
(\ref{fuchs1})
of the same form (\ref{local}) analytic on
\beq\label{cuts1}
{\mathbb C}\setminus \cup_j \hat L_j.
\eeq
Here we choose the branches of the functions $(u_j-\lambda)^{\nu-{1\over 2}}$
near $\hat L_j$ by rotating the branches of the same functions defined near
$L_j$ in the counter-clockwise direction.

\begin{lemma} \label{lemma36} The result of the counter-clockwise 
analytic continuation of the solutions
$\phi^{(1)}$, \dots, $\phi^{(n)}$ is related to the solutions $\hat\phi^{(1)}$,
\dots, $\hat\phi^{(n)}$ by the transformation
\beq\label{left-right}
\left( \phi^{(1)}, \dots, \phi^{(n)}\right) =
\left( \hat\phi^{(1)}, \dots, \hat\phi^{(n)}\right) 
\left[ 1+(1+q ^{-1}) G_+\right]
\eeq
Similarly, for the clockwise analytic continuation of the same functions
one obtains
\beq\label{clock}
\left( \phi^{(1)}, \dots, \phi^{(n)}\right) = -q^{-1}
\left( \hat\phi^{(1)}, \dots, \hat\phi^{(n)}\right) 
\left[ 1+(1+q ) G_-\right].
\eeq
Here the matrices $G_\pm$ are defined as follows
\eqa\label{g-pm}
&&
\left( G_+\right)^{ij}=\left\{ \begin{array}{ll} 0,  & i\geq j \\
g^{ij}, 
& i<j \end{array}\right.
\nn\\
&&
\left( G_-\right)^{ij} = \left\{ \begin{array}{ll} 0, & i\leq j \\
g^{ij}, & i>j \end{array}\right.
\\
\eeqa
\end{lemma}

\pf Let us prove the first formula. Counter-clockwise 
analytic continuation of $\phi^{(1)}$ till $\hat L_1$ 
does not meet obstructions: we simply rotate the branchcut $L_1$ untill $\hat
L_1$. So
$$
\hat\phi^{(1)}=\phi^{(1)}.
$$
To continue analytically $\phi^{(2)}$ till $\hat L_2$
we are to cross the branchcut $L_1$. This changes $\phi^{(2)}$ to
$M_1 ^{-1}\phi^{(2)}$. Using (\ref{anti-mono}) for $i=2$, $j=1$ we obtain
for the analytic continuation
$$
\hat\phi^{(2)}=M_1^{-1} \phi^{(2)}=\phi^{(2)}-(1+q^{-1}) g^{12} \hat\phi^{(1)}.
$$
Similarly, 
$$
\hat\phi^{(3)} =M_1^{-1}M_2^{-1} \phi^{(3)} =
M_1^{-1} \phi^{(3)} -(1+q^{-1}) g^{23} M_1^{-1}\phi^{(2)}
$$
$$
=\phi^{(3)} -(1+q^{-1})\left[ \hat\phi^{(1)} g^{13} +g^{23} \hat\phi^{(2)}\right],
$$
etc.,
$$
\hat\phi^{(n)} = M_1^{-1}M_2^{-1}\dots M_{n-1}^{-1}\phi^{(n)} =
\phi^{(n)} -(1+q^{-1}) \left[ \hat\phi^{(1)}g^{1n} +\hat\phi^{(2)}g^{2n} +\dots
+ \hat\phi^{(n-1)} g^{n-1,n}\right].
$$
This gives (\ref{left-right}). Similarly, continuing analytically the functions $\phi^{(i)}$
in the clockwise direction and using (\ref{monoij}) we derive (\ref{clock})
(here we are to take into account the change of the branches of the functions
$(u_i-\lambda)^{\nu-{1\over 2}}$). Lemma is proved.

Let us define now vector-functions $Y_R^{(i)}(z)$ and $Y_L^{(i)}(z)$ via the
following (inverse) Laplace transforms
\eqa\label{laplace}
&&
Y_R^{(i)}(z) := 
{i\over (2\pi)^{3/2}}{q\over 1+q} z^{\nu+{1\over 2}}\oint_{\hat L_i} \hat\phi^{(i)}(\lambda)e^{z\,\lambda}d\lambda
\nn\\
&&
Y_L^{(i)}(z) :={i\over (2\pi)^{3/2}}{q\over 1+q} z^{\nu+{1\over 2}}
\oint_{L_i} \phi^{(i)}(\lambda)e^{z\,\lambda}d\lambda.
\eeqa
Here the loop integrals are taken along the infinite cycle coming from infinity
along the left shore of the branchcut $L_j / \hat L_j$ resp., then encircling
the point $\lambda=u_j$ and returning to infinity along another shore of the
same branchcut. We define the branch of $z^{\nu +{1\over 2}}$ doing a branchcut
along the negative half of the line $\ell$.

\begin{lemma} \label{lemma37} The vector-functions $Y_{R/L}^{(i)}(z)$ are linearly
independent
solutions to the 
system (\ref{isomono}). They are analytic in the half-planes $\Pi_{R/L}$ resp. to the
right/to the left of the line $\ell$. i.e.,
\beq\label{pi}
\Pi_+ =\{ z~ | ~ \varphi < \arg \,z< \varphi+\pi\}.
\eeq
In these half-planes they have
the asymptotic development of the form
\beq\label{asy}
Y_{R/L}^{(i)}(z) \sim \left( e_i +O\left( 1/z\right)\right) e^{z\, u_i}.
\eeq
\end{lemma}

\begin{cor} The Stokes matrix of the operator (\ref{oper}) w.r.t. the 
oriented line
$\ell$ is equal to
\beq\label{stokes}
S= 1+(1+q^{-1}) G_+.
\eeq
The transposed Stokes matrix is equal to
\beq\label{stokes-t}
S^T =1+(1+q)G_-.
\eeq
\end{cor}

Proof of Lemma. At infinity the Fuchsian system (\ref{fuchs1}) has a regular singularity. 
So the solutions
$\phi^{(i)}(\lambda)$ and $\hat\phi^{(i)}(\lambda)$ grow at infinity
not faster than some power of $|\lambda|$. This proves convergence of the
integrals for $|\lambda|\to \infty$, $z\in \Pi_{R/L}$ resp. 
Analyticity of the integrals (\ref{laplace}) in the half-planes
$\Pi_R$/$\Pi_L$ resp. is a standard fact of the theory of Laplace integrals
(see, e.g., \cite{lap}). Integrating the convergent expansions near $\lambda=u_i$
of the form (\ref{local})
of the functions
$\phi^{(i)}(\lambda)$/$\hat\phi^{(i)}(\lambda)$ multiplied by $e^{z\, \lambda}$
we 
arrive at the asymptotic developments (\ref{asy}) at $z=\infty$ of 
the integrals.
Plugging the integrals (\ref{laplace}) into (\ref{isomono}) and integrating 
by parts (we can neglect
the boundary terms due to the exponential vanishing of the integrand
at infinity) we prove that the integrals satisfy the system (\ref{isomono}).
Their linear independence follows from the independence of the principal terms
of the asymptotic developments. Lemma is proved.\epf

Proof of Corollary. To analytically continue the integral $Y_L^{(i)}(z)$
in the clockwise direction to the half-plane $\Pi_L$ through
the positive part of the line $\ell$ one is to rotate 
counter-clockwise the contour
$L_i$ till it will get to the position $\hat L_i$. Using (\ref{left-right}) we derive
that in a narrow sector around the positive half-line $\ell$
$$
\left(Y_L^{(1)}, \dots, Y_L^{(n)}\right) = 
\left(Y_R^{(1)}, \dots, Y_R^{(n)}\right)\left[ 1+(1+q^{-1}) G_+\right] .
$$
This gives the formula (\ref{stokes}) for the Stokes matrix. 
Similarly, continuing analytically $Y_L^{(i)}(z)$ in the counter-clockwise
direction through the negative part of the line $\ell$ (here the branch
of $z^{\nu+{1\over 2}}$ changes) and using (\ref{clock}) we arrive at 
$$
\left(Y_L^{(1)}, \dots, Y_L^{(n)}\right) = 
\left(Y_R^{(1)}, \dots, Y_R^{(n)}\right)\left[ 1+(1+q) G_-\right] .
$$
As we know from the theory of Stokes matrices for the operator (\ref{oper}), in the narrow
sector near the positive part of $\ell$ 
$$
\left(Y_L^{(1)}, \dots, Y_L^{(n)}\right) = 
\left(Y_R^{(1)}, \dots, Y_R^{(n)}\right) S
$$
and in the narrow sector near the negative part of $\ell$
$$
\left(Y_L^{(1)}, \dots, Y_L^{(n)}\right) = 
\left(Y_R^{(1)}, \dots, Y_R^{(n)}\right) S^T
$$
(see, e.g., \cite{pain}). Corollary is proved.\epf

Using (\ref{stokes}) and (\ref{stokes-t}) we rewrite the formula 
(\ref{monoij}) of the monodromy transformation that we redenote $M_i(q)\in
GL(n, {\mathbb C}[q, q^{-1}])$
as follows
\beq\label{mono-stokes}
M_i(q)\phi^{(j)} =\left\{ \begin{array}{lll}
\phi^{(j)} - q\,s_{ij} \phi^{(i)}, & i<j\\
-q\, \phi^{(i)}, & i=j \\
\phi^{(j)} - s_{ji} \phi^{(i)},  & i>j
\end{array}\right.
\eeq
This is the reflection (\ref{reflect}) w.r.t. the bilinear form (\ref{q-form}). 

It remains to prove linear independence of the solutions $\phi^{(1)}$, \dots,
$\phi^{(n)}$ under the assumption (\ref{half}). Any linear dependence
$$
c_1 \phi^{(1)} +\dots + c_n\phi^{(n)}=0
$$
would give a vector invariant w.r.t. to the total monodromy operator
$M_n M_{n-1}\dots M_1$ corresponding to a big counter-clockwise loop
around the origin. From (\ref{left-right}), (\ref{clock}) it follows that such an operator 
acts as follows
\beq\label{total}
\left( \phi^{(1)}, \dots, \phi^{(n)}\right) \mapsto -q
\left( \phi^{(1)}, \dots, \phi^{(n)}\right)
S^{-T} S.
\eeq
Here we denote 
$$
S^{-T}:= \left( S^T\right)^{-1}.
$$
The transformation (\ref{total}) has an invariant vector {\it iff}
$$
\det \left[q\, S^{-T} S+1\right] =0.
$$
This coincides with degeneracy of the matrix (\ref{q-form}). The contradiction proves
independence of the solutions.

To complete the proof of Theorem in the non-resonant case we are to prove
that the solutions $\phi^{(i)}$ satisfy also the equations (\ref{fuchs3}). 
To this end we consider the fundamental matrix 
\beq\label{fund}
\Phi(\lambda;u):=
\left( \phi^{(1)}(\lambda; u), \dots, \phi^{(n)}(\lambda;u)\right)
\eeq
(recall that the dependence of this matrix on $u=(u_1, \dots, u_n)$ is
determined uniquely). Due to compatibility of the system (\ref{fuchs1}),
(\ref{fuchs3}) the matrix
$$
\partial_i\Phi -\left( {E_i \left( {1\over 2} -\nu +V\right)\over \lambda -u_i}
+V_i\right) \Phi
$$
is again a matrix solution of (\ref{fuchs1}). So it has the form
$$
\partial_i\Phi -\left( {E_i \left( {1\over 2} -\nu +V\right)\over \lambda -u_i}
+V_i\right) \Phi=\Phi\, T_i
$$
for some matrix $T_i$ independent on $\lambda$. Using the expansions
$$
\phi^{(i)} =(u_i-\lambda)^{\nu-{1\over 2}}
\left[ e_i +{e_i\, (u_i-\lambda)\over \nu+{1\over 2}} (V_i e_i - 2 H_i e_i)
+O\left( u_i-\lambda\right)^2\right], ~~\lambda\to u_i
$$
with
$$
H_i={1\over 2} \sum_{j\neq i} {V_{ij}^2\over u_i-u_j}
$$
of the solutions (\ref{local}) and
$$
\psi = a+O(u_i-\lambda), ~~\lambda\to u_i, ~~\left( \nu-{1\over 2}+V\right)\,a=0
$$
of any solution to (\ref{fuchs1}) analytic at $\lambda=u_i$ we prove that $T_i=0$.
This completes the proof of Theorem in the nonresonant case (\ref{half}).

Before going further we will give an interpretation of the bilinear form 
$(~,~)_\nu$ defined in (\ref{q-form}).
Let us denote $L(\nu)$ the $n$-dimensional space of solutions
to the system (\ref{fuchs1}), (\ref{fuchs3}). We have a natural pairing (cf.
(\ref{pairing}))
\eqa\label{pairing1}
&&
L(-\nu)\times L(\nu)\to {\mathbb C},
\nn\\
&&
\left( \phi_{-\nu}, \psi_{\nu} \right) := \phi_{-\nu}^T(\lambda; u) 
(U-\lambda)
\psi_{\nu}(\lambda;u), ~~\psi_\nu\in L(\nu), ~\phi_{-\nu}\in L(-\nu).
\eeqa
It is easy to see that the pairing does not depend on $\lambda$ neither on $u$.
Clearly, this pairing does not degenerate under the nonresonancy assumption
(\ref{half}).

\begin{lemma} In the nonresonant case the matrix of the pairing w.r.t.
the bases $\phi^{(1)}_\nu$, \dots, $\phi^{(n)}_\nu$ and
$\phi^{(1)}_{-\nu}$, \dots, $\phi^{(n)}_{-\nu}$ in $L(\nu)$ and $L(-\nu)$ resp.
coincides, up to a scalar factor with (\ref{q-form}):
\beq\label{q-pair}
\left( \phi_{-\nu}^{(i)}, \phi_{\nu}^{(j)}\right) =\left( q^{1/2} S +
q^{-1/2}S^T\right)_{ij}.
\eeq
So, the formula for the monodromy transformations (\ref{monoij}) can be recast
into the form
\beq\label{mono-q}
M_i \phi_\nu^{(j)} =\phi_\nu^{(j)} -q^{1/2}\left( \phi_{-\nu}^{(i)},
\phi_\nu^{(j)}\right) \, \phi_\nu^{(i)}.
\eeq
\end{lemma}

\pf We observe first that, for any $i$ the solution $\phi_\nu^{(i)}$
is orthogonal to any solution $\psi\in L(-\nu)$  analytic at
$\lambda=u_i$. So, using (\ref{gij}) we obtain
$$
\left( \phi_{-\nu}^{(i)}, \phi_{\nu}^{(j)}\right)=q^{-1/2}(1+q)
g^{ij}.
$$
Using (\ref{stokes}), (\ref{stokes-t}) we arrive at the proof of Lemma. \epf

Let us now proceed to the resonant case. We consider first the particular case
$\nu={1\over 2}$ (actually, this is still a nonresonant case in the standard
sense of the theory of linear ODEs with rational coefficients. However,
the monodromy matrices appear to have nontrivial Jordan blocks). The matrices
$A_i$ are now all nilpotent,
$$
A_i^2=0.
$$
So they have just one block $\left(\begin{matrix} 0 & 1 \cr 0 & 0 \cr
\end{matrix}\right)$
in their Jordan normal form. The monodromy matrix $M_i$ has $n-1$ linearly 
independent eigenvectors with the eigenvalue 1 and, if the $i$-th row of the
matrix $V$ is nonzero, it has
one Jordan block
$\left(\begin{matrix} 1 & 2\pi i \cr 0 & 1 \cr\end{matrix}\right)$. That means that there are
$n-1$ solutions to (\ref{fuchs1}) analytic at $\lambda=u_i$. One of these solutions
belongs to the image of $M_i -{\rm id}$. We denote it, like above, $\phi^{(i)}$.
It can be normalized in such a way that
\beq\label{phii}
\phi^{(i)}(\lambda) =e_i +O(u_i-\lambda), ~~\lambda\to u_i.
\eeq
Such a normalization determines the solution uniquely.
A logarithmic solution denoted by $\chi^{(i)}(\lambda)$ 
corresponds to $\phi^{(i)}(\lambda)$ in the following sense
\beq\label{mi}
\left(M_i-{\rm id}\right) \chi^{(i)} =2\pi\,i\,\phi^{(i)}.
\eeq
Near $\lambda=u_i$ it can be written as
\beq\label{chi}
\chi^{(i)}(\lambda) =\log (u_i-\lambda) \,\phi^{(i)}(\lambda)
+\delta^{(i)}(\lambda)
\eeq
where $\delta^{(i)}(\lambda)$ is a vector-function of $\lambda$ analytic at 
$\lambda=u_i$. It is easy to see that the value of this function at
$\lambda=u_i$ must satisfy the condition
\beq\label{cond-i}
\left[ V\delta^{(i)}(\lambda=u_i)\right]_i=-1
\eeq
(here and below $[~]_i$ stands for the $i$-th coordinate of the vector). 
To obtain a basis in the space  of solutions to (\ref{fuchs1}) 
with $\nu=1/2$ we are to add
to $\phi^{(i)}$ and $\chi^{(i)}$, where $i$ is fixed, $n-2$ independent 
solutions analytic
at $\lambda=u_i$. Any such a solution $\psi(\lambda)$ must satisfy, like in
(\ref{cond-i}),
the equation 
$$
\left[ V\psi(\lambda=u_i)\right]_i=0.
$$
We can always assume that the $i$-th coordinate of $\psi(u_i) $ is equal
to zero.

Globally the solutions $\phi^{(1)}$, \dots, $\phi^{(n)}$ are analytic in
(\ref{cuts})
with the same choice of the branch cuts $L_1$, \dots, $L_n$ as above.
We define a matrix $g^{ij}$ in such a way that $g^{ii}=0$ and
\beq\label{phi-chi}
\phi^{(i)}(\lambda)={1\over 2\pi\,i} g^{ji}\chi^{(j)}(\lambda) +\, analytic,
~~\lambda\to u_j, ~~i\neq j.
\eeq
Like in Lemma \ref{lemma35}, using the last equation we prove

\begin{lemma} The monodromy transformation $M_i$ acts onto the solution
$\phi^{(j)}$ as folllows
\beq\label{mij}
M_i\phi^{(j)} =\phi^{(j)} + g^{ij} \phi^{(i)}.
\eeq
\end{lemma}

Observe that
$$
M_i^{-1}\phi^{(j)} =\phi^{(j)} - g^{ij} \phi^{(i)}.
$$

Let us introduce another system $\hat\phi^{(1)}$, \dots, $\hat\phi^{(n)}$
of solutions to (\ref{fuchs1}) with the same behaviour (\ref{phii})
but with the branchcuts along $\hat L_1$, \dots, $\hat L_n$.

\begin{lemma} The results of counter-clockwise/clockwise analytic
continuation of the solutions $\phi^{(1)}$, \dots, $\phi^{(n)}$ are given by the
following formulae
\beq\label{gpm}
\left(\phi^{(1)}, \dots, \phi^{(n)}\right) =\left(
\hat\phi^{(1)}, \dots, \hat\phi^{(n)}\right) (1\pm G_\pm).
\eeq
Here $G_+$/$G_-$ are the upper/lower triangular parts of the matrix
$G=\left(g^{ij}\right)$.
\end{lemma}

Proof is similar to that of Lemma \ref{lemma36}.

We introduce also the solutions $\hat\chi^{(1)}$, \dots, $\hat\chi^{(n)}$
related to $\hat\phi^{(1)}$, \dots, $\hat\phi^{(n)}$ by the formulae 
(\ref{chi}).

\begin{lemma} \label{lemma312} The matrix $G$ has the form
\beq\label{anti-g}
G=S-S^T.
\eeq
\end{lemma}

\pf We construct solutions to the system (\ref{isomono}) using the following Laplace
integrals
\eqa\label{laplace1}
&&
Y_L^{(i)}(z) ={z\over 2\pi\,i} \oint _{L_i}
\chi^{(i)}(\lambda)e^{\lambda\,z} d\,\lambda = -z\int_{L_i} \phi^{(i)}(\lambda)
e^{\lambda\,z}d\,\lambda
\nn\\
&&
Y_R^{(i)}(z) ={z\over 2\pi\,i} \oint _{\hat L_i}
\hat\chi^{(i)}(\lambda)e^{\lambda\,z} d\,\lambda
=-z \int_{\hat L_i}\phi^{(i)}(\lambda).
\eeqa
Like in Lemma \ref{lemma37} we prove that these are the columns of the matrix solutions
$Y_{L/R}(z)$ to (\ref{isomono}) having the asymptotic development (\ref{asy}) 
in $\Pi_{L/R}$ resp.
From (\ref{gpm}) we conclude that 
$$
S=1+G_+, ~~S^T=1-G_-.
$$
This proves Lemma.\epf

We define now an antisymmetric bilinear form on $L(1/2)$ by one of the following
three expressions:
\beq\label{poi}
\left( \phi, \psi\right)_{1/2} := -2\pi\, \phi^T(\lambda) \, V \, \psi(\lambda)
= 2\pi\, {d\phi^T\over d\lambda} \left(U-\lambda\right) \psi(\lambda)
= -2\pi\, \phi^T(\lambda) \left(U-\lambda\right) {d\psi\over d\lambda}.
\eeq
It is induced by the pairing $L\left( -{1\over 2}\right) \times
 L\left( {1\over
2}\right) \to {\mathbb C}$ and by the shift operator $L\left( -{1\over 2}\right) 
\to L\left( {1\over
2}\right) $ defined by (\ref{shift}).

\begin{lemma} The matrix of the bilinear form is
\beq\label{mat-poi}
\left( \phi^{(i)}, \phi^{(j)}\right)_{1/2} =-i\left( S-S^T\right)_{ij}.
\eeq
\end{lemma}

\pf It is easy to see that $\left(\phi^{(i)}, \psi\right)_{1/2}=0$
for any solution $\psi(\lambda)$ analytic at $\lambda=u_i$. Using (\ref{cond-i})
we obtain that $\left(\phi^{(i)}, \chi^{(i)}\right)_{1/2}=2\pi$. 
From these two facts and from (\ref{phi-chi}) we derive that
$$
\left(\phi^{(i)}, \phi^{(j)}\right)_{1/2}=-i \, g^{ij}.
$$
From Lemma \ref{lemma312} we obtain (\ref{mat-poi}). Lemma is proved.\epf

As a consequence we obtain that, under the assumption (\ref{ch-q}), (i.e., the
assumption of nondegeneracy of the matrix $S-S^T$) 
the solutions $\phi^{(1)}$, \dots,  $\phi^{(n)}$ form a basis in the space 
$L(1/2)$,
and that the monodromy transformations act in this basis as in (\ref{reflect}):
$$
M_i\phi^{(j)} =\phi^{(j)}-i\left(\phi^{(i)},\phi^{(j)}\right)_{1/2} \,\phi^{(i)}.
$$
As above we prove that the functions $\phi^{(1)}(\lambda;u)$, \dots,  
$\phi^{(n)}(\lambda;u)$ satisfy also the equations (\ref{fuchs3}).

We have proved Theorem for $\nu=1/2$. To obtain the proof for $\nu=-m+{1\over
2}$ with a positive integer $m$ we just use the isomorphism
$$
{d^m\over d\lambda^m}\,: \, L(1/2) \to L(-m+1/2).
$$
This produces a basis in $L(-m+1/2)$ with the needed monodromy. Finally,
for $\nu=m+{1\over 2}$ with a positive integer $m$ we use the non-degenerate
pairing (\ref{pairing1}) to construct a basis in $L(m+1/2)$
dual to the basis in $L(-m+1/2)$. This completes the proof of Theorem.\epf

We will now use the above Theorem in order to describe the geometry of the
discriminant $\Sigma$ w.r.t. the geometry induced by the intersection form.
We first recall that, in any chart $Fr(\hat\mu, R, e, S, C)$ on $M_{s\, s}$ the
intersection $\Sigma\cap M_{s\, s}$ splits into the union of hypersurfaces
\beq\label{uk=0}
\Sigma\cap M_{s\, s} = \cup_{k=1}^n \left\{ u_k =0 \right\} , \quad u_i\neq u_j
\quad {\rm for} ~i\neq j.
\eeq
We will describe the behaviour of the flat coordinates $p^1$, \dots, $p^n$ of the intersection form
on (\ref{uk=0}). Within any connected and simply connected domain in
the coordinate patch $Fr(\hat\mu, R, e, S, C)\setminus \Sigma$ one can choose, for
$d\neq 1$ the
particular system of coordinates by putting
\beq\label{p-choice}
p^a (u)= {2\over 1-d} \sum_{i=1}^n u_i \psi_{i1}(u) \phi_i^{(a)}
(u; \lambda=0;\nu=0), \quad a=1, \dots, n
\eeq
(cf. the formula (\ref{reco})). Here the basis 
$\phi^{(a)}= (\phi_1^{(a)}, \dots,\phi^{(a)}_n)^T$ of solutions to the 
Fuchsian system (\ref{fuchs1}) is chosen as
in (\ref{local}). Note that, due to (\ref{q-pair}) (for $\nu=0$) one has
\beq\label{gram0}
(dp^a, dp^b) =:G^{ab} = \left( S+S^T\right)^{ab}.
\eeq

We will describe the limiting behaviour of the basis of the
periods (\ref{p-choice}) near (\ref{uk=0}).
The following statement is a refinement of Lemma G.2 in \cite{cime}.

\begin{theorem} Let $M$ be a semisimple Frobenius manifold satisfying
$d\neq 1$ and
\beq\label{euc}
\det (S+S^T)\neq 0.
\eeq
Let ${\bf D}\subset Fr(\hat\mu, R, e, S, C)\setminus \Sigma\subset M_{s\,
s}\setminus \Sigma$
be a connected simply connected domain such that $\Sigma\subset \overline {\bf
D}$. Then the functions $p^1(u)$, \dots, $p^n(u)$ can be extended continuously
up to $\Sigma$. With respect to this continuation the component
$$
u_k=0
$$
of $\Sigma$ becomes an affine hyperplane
\beq\label{pk=0}
p^k=0.
\eeq
\end{theorem}

\pf It is technically more convenient to compute the limiting behaviour of
the $\lambda$-periods
\beq\label{p-choice1}
\tilde p^a(u;\lambda)
={2\over 1-d} \sum_{i=1}^n (u_i-\lambda) \psi_{i1}(u) \phi_i^{(a)}
(u; \lambda;\nu=0), \quad a=1, \dots, n
\eeq
on $\Sigma_\lambda$. Because of (\ref{gram0}), \ref{euc}) the functions (\ref{p-choice1})
are independent.
Note that the $k$-th component of the intersection of 
$\Sigma_\lambda\cap M_{s\,s}$ are the hypersurfaces 
$$
u_k=\lambda.
$$
Near $u_k=\lambda$ one has
\beq\label{beh1}
\phi^{(k)}_i = {\sqrt{2}\over \sqrt{u_k-\lambda}}\left( \delta^k_i +
O(u_k-\lambda)\right)
\eeq
and, for $l\neq k$
\beq\label{beh2}
\phi^{(k)}_i= {G^{kl}\over \sqrt{u_l-\lambda}} \left( \delta^l_i +
O(\sqrt{u_l-\lambda})\right).
\eeq
So, near $u_k=\lambda$
\beq\label{p-beh1}
\tilde p^k(u;\lambda) = {2\sqrt{2}\over 1-d} \psi_{k1}(u) \sqrt{u_k-\lambda}
+O(u_k-\lambda)
\eeq
and, for $l\neq k$
\beq\label{p-beh2}
\tilde p^l(u;\lambda) =\tilde p^l_0 + {\sqrt{2}\over 1-d} G^{kl} \psi_{k1}(u)
\sqrt{u_k-\lambda} + O(u_k-\lambda)
\eeq
where $p^l_0=p^l_0(u_1, \dots, \hat u_k, \dots, u_n; \lambda)$ is an analytic
function on 
$$
Fr(\hat\mu, R, e, S, C)\setminus \left\{ u_k=\lambda \right\} .
$$
From these formulae the continuity of the functions (\ref{p-choice1}) up to 
$\Sigma_\lambda$ readily follows. In particular, for $\lambda=0$ one obtains
continuity of the periods (\ref{p-choice}) up to $\Sigma$. Also (\ref{pk=0})
readily follows from (\ref{p-beh1}). The Theorem is proved.\epf

Observe that the angles between the hyperplanes (\ref{pk=0}) can be computed
from the Gram matrix (\ref{gram0}).

We will now apply the formulae (\ref{p-beh1}), (\ref{p-beh2}) in order to
describe the analytic properties of the dual almost Frobenius manifold
near $\Sigma$.

\begin{theorem}\label{c-th} Under the assumptions $d\neq 1$ and (\ref{euc}) the structure
coefficients $\cstar_{ab}^c(p)$ of the dual Frobenius manifold are continuous
up to the hyperplane (\ref{pk=0}) for all $a$, $b$, $c$ except
\beq\label{c-beh}
\cstar_{kk}^a = {1\over 1-d} {1\over p^k} G^{ka} +\, {\rm regular} ~ {\rm
terms}, \quad p^k\to 0.
\eeq
\end{theorem}

\pf Let us first derive a formula for $\cstar_{ab}^c(u)$. We will use the
coordinates (\ref{p-choice1}) and we will set $\lambda=0$ at the end of the
computation.

\begin{lemma} At any point $u=(u_1, \dots, u_n)\in M_{s\, s}\setminus
\Sigma_\lambda$ one has
$$
{\pal\over \pal p^b} * {\pal\over \pal p^c} = \cstar_{bc}^a(u) {\pal\over \pal
p^c}
$$
where
\beq\label{c-formula}
\cstar_{bc}^a(u)= G_{cd}G_{bf} \sum_{i=1}^n {u_i-\lambda\over\psi_{i1}(u)}
\phi^{(a)}_i(u;\lambda)\phi^{(f)}_i(u;\lambda)\phi^{(d)}_i(u;\lambda).
\eeq
\end{lemma}

\pf From (\ref{p-phi}) it follows that
$$
d\tilde p^a = \sum_{i=1}^n \psi_{i1}(u) \phi_i^{(a)}(u;\lambda) du_i.
$$
From the definition of the canonical coordinates it follows that
$$
du_i \cdot du_j =\delta_{ij}\psi_{i1}^{-2} du_i.
$$
So
$$
d\tilde p^a \cdot d\tilde p^b = \sum_{i=1}^n  \phi_i^{(a)} \phi_i^{(b)} du_i
$$
(we will omit writing the arguments $u$, $\lambda$ of the functions under 
consideration for the sake of brevity of the formulae). Using
$$
\sum_{i=1}^n (u_i-\lambda) \phi_i^{(a)} \phi_i^{(b)} = G^{ab}
$$
we derive that
$$
du_i = \psi_{i1}^{-1} (u_i-\lambda) G_{ab} \phi_i^{(a)} dp^b.
$$
So
$$
d\tilde p^a\cdot d\tilde p^b = G_{cd}\sum_{i=1}^n \psi_{i1}^{-1} (u_i-\lambda)
\phi_i^{(a)} \phi_i^{(b)} \phi_i^{(d)}\, d\tilde p^c.
$$
In other words, in the coordinates (\ref{p-choice1}) the structure coefficients
of the multiplication $\cdot$ on the cotangent bundle read
$$
c^{ab}_c (u) = G_{cd}\sum_{i=1}^n \psi_{i1}^{-1} (u_i-\lambda)
\phi_i^{(a)} \phi_i^{(b)} \phi_i^{(d)}.
$$
But we already now that on the cotangent bundle 
$$
\cstar^{ab}_c = c^{ab}_c.
$$
Lowering the index
$$
\cstar_{bc}^a = G_{bf}\cstar^{fa}_c
$$
we obtain the needed formula. The Lemma is proved. \epf

Prove of the Theorem. Substituting the expansions (\ref{p-beh1}), (\ref{p-beh2})
into the formula (\ref{c-formula}) we obtain, near $u_k=\lambda$
$$
\cstar_{bc}^a = G_{cd}G_{bf} \sum_{i=1}^n {u_i-\lambda\over 2\sqrt{2} \psi_{i1}
(u_k-\lambda)^{3/2}} G^{ak} G^{fk} G^{dk} \left[ \delta_i^k
+O(\sqrt{u_k-\lambda})\right]
$$
$$
= \delta_c^k \delta_b^k {G^{ak}\over 2\sqrt{2} \psi_{k1} \sqrt{u_k-\lambda}}
+O(1).
$$
From this formula and from (\ref{p-beh1}) the asymptotic formula 
(\ref{c-beh}) immediately follows.\epf

Similarly to (\ref{p-choice}) one can introduce, assuming $d\neq 1$ a system of
$n$ odd periods 
\beq\label{pi-choice}
\varpi^a ={2\over 1-d} \sum_{i=1}^n u_i \psi_{i1}(u) \phi_i^{(a)}
(u; \lambda=0;\nu={1\over 2}), \quad a=1, \dots, n.
\eeq
From (\ref{mat-poi}) it follows

\begin{cor} The functions $\varpi^1$, \dots, $\varpi$ have the following
constant matrix of the Poisson brackets (\ref{pb})
\beq\label{pi-pb}
<d\varpi^a {\mathcal V}, d\varpi^b> = i \, (S-S^T)_{ab}.
\eeq
\end{cor}

We will finally describe the behaviour of the fundamental system $\left( \phi^{(1)},
\dots, \phi^{(n)}\right)$ at $|\lambda|\to\infty$. For simplicity only the
generic case
(\ref{half}) will be considered.

Before we proceed to the formulation of the Theorem, let us recall some useful
formulae from the theory of Laplace integrals. First, 
\beq\label{ga0}
\int_0^\infty z^{s-{1\over 2}} e^{-\lambda\, z} dz = \Gamma\left( s+{1\over
2}\right) \lambda^{-s-{1\over 2}}.
\eeq
This formula coincides with the definition of the gamma function for ${\rm Re}\,
s > -{1\over 2}$; for other $s\ni {\mathbb Z}+{1\over 2}$ it is obtained by
analytic continuation. Differentiating, we derive
\beq\label{ga1}
\int_0^\infty z^{s-{1\over 2}} \log^k ze^{-\lambda\, z} dz = {d^k\over ds^k}\Gamma\left( s+{1\over
2}\right) \lambda^{-s-{1\over 2}}.
\eeq
We will also need a matrix analogue of the last formula.

\begin{lemma} The following formula holds true
\beq\label{matrix}
\int_0^{\infty}e^{-\lambda\, z} z^{\hat\mu + s -{1\over 2}}
z^R dz = \sum_{m \geq 0} \left[e^{R\, \pal_s}\right]_m\Gamma \left( \hat\mu + s+m +{1\over 2}\right)
\times \lambda^{-(s+m+\hat\mu+{1\over 2})}
\lambda^{-R}. 
\eeq
\end{lemma}

In this formula the $m$-th component $\left[ ~~\right]_m$
of the matrix is defined like in (\ref{comp}).

\pf We have
$$
\int_0^{\infty }e^{-\lambda\, z} z^{s+ \hat\mu  -{1\over 2}}
z^R dz = \int_0^\infty e^{-t} t^{\hat\mu + s -{1\over 2}}
t^{R_0 +{R_1\over \lambda} +{R_2\over \lambda^2}+\dots} dt \, 
\lambda^{-(\hat\mu
+s+{1\over 2})} \lambda^{-R}
$$
$$
=\sum_{m\geq 0}\sum_{k\geq 0} \int_0^\infty e^{-t} t^{\hat\mu + s-{1\over 2}}
\lambda^{-m}{[R^k]_m \log^k t\over k! } dt \, \lambda^{-(\hat\mu
+s +{1\over 2})} \lambda^{-R}
$$
$$
=\sum_{m\geq 0}\sum_{k\geq 0} {1\over k!} \pal_s^k \Gamma(\hat\mu+s+{1\over
2}) {[R^k]_m\over \lambda^m}\lambda^{-(\hat\mu
+s +{1\over 2})} \lambda^{-R}
$$
$$
= \sum_{m\geq 0} \sum_{k\geq 0} \left[ {R^k\over k!}\right]_m
\Gamma(\hat\mu+s+m+{1\over 2}) \times \lambda^{-(\hat\mu
+s+m +{1\over 2})}\lambda^{-R}
$$
$$
= \sum_m 
[e^{R \partial_s}]_m\Gamma(\hat\mu+s+m+{1\over 2}) \times \lambda^{-(\hat\mu+
s+m+{1\over 2})} \lambda^{-R}. 
$$
The Lemma is proved. \epf

\begin{theorem} At $|\lambda|\to\infty$ within the domain (\ref{cuts}) the
fundamental system of solutions (\ref{local}) has the following expansion
\eqa\label{inf}
&&
\Phi= \left( \phi^{(1)},
\dots, \phi^{(n)}\right) 
\nn\\
&&
={i\over \sqrt{2\pi}}(q^{1/2}+q^{-1/2})^{-1} 
\sum_{p=0}^\infty\sum_{m\geq 0} \Theta_p(t) \left[ e^{-R\,\pal_\nu}\right]_m
\Gamma \left(
p+m+\hat\mu-\nu +{1\over 2}\right)
\nn\\
&&
\times\lambda^{-(p+m+\hat\mu-\nu+{1\over 2})}
\lambda^{-R} C^{-1}\,\left( q^{1/2}S+q^{-1/2}S^T\right).
\nn\\
\eeqa
\end{theorem}

\pf For ${\rm Re}\,\nu <<0$ let us consider the matrix
$$
\tilde\Phi(t;\lambda) ={i\over \sqrt{2\pi}} (1+q^{-1}) 
\int_0^\infty Y_R(z;t) e^{-\lambda\, z} {dz\over z^{\nu+{1\over 2}}}.
$$
Here the integration is to be performed along a ray on the $z$-plane belonging
to $\Pi_R$. It is a fundamental matrix of solutions to (\ref{fuchs1}) analytic in
(\ref{cuts}) satisfying the following property: the $i$-th column of $\Phi(t;\lambda)$
near $\lambda=u_i$ behaves
$$
\left( \Phi_0(t;\lambda)\right)_i = \phi^{(i)}(t;\lambda) + analytic
$$
and it is analytic near $\lambda=u_j$ for any $j\neq i$. Using (\ref{gij}),
(\ref{stokes}), (\ref{stokes-t}) we
obtain
\beq\label{til-phi}
\Phi_0(t;\lambda) =(1+q) \left( \phi^{(1)}(t;\lambda), \dots,
\phi^{(n)}(t;\lambda)\right) \left( q\, S+S^T\right)^{-1}.
\eeq
On the other hand, we can obtain the expansion of $\tilde\Phi(t;\lambda)$ near the
regular singularity
$\lambda=\infty$ replacing $Y_L$ by $Y_0C$ and integrating the
series (\ref{theta1})
\eqa\label{riad}
&&
{i\over \sqrt{2\pi}} (1+q^{-1}) 
\int_0^\infty Y_0(z;t) e^{-\lambda\, z} {dz\over z^{\nu+{1\over 2}}}
\nn\\
&&
={i\over \sqrt{2\pi}} (1+q^{-1}) 
\int_0^\infty \sum_p \Theta_p(t) z^{p+\hat\mu-\nu-{1\over 2}} z^R e^{-\lambda\, z} 
dz
\nn\\
&&
={i\over \sqrt{2\pi}}(1+q^{-1}) \sum_p \sum_{m\geq 0}\Theta_p(t) \left[ e^{-R\,
\pal_\nu}\right]_m\Gamma(p+m+\hat\mu-\nu+{1\over 2})
\nn\\
&& 
\times \lambda^{-(p+m+\hat\mu-\nu+{1\over 2} } \lambda^{-R}.
\eeqa
Comparing (\ref{til-phi}) and (\ref{riad}) we obtain (\ref{riad}). 
Theorem is proved.\epf

We will now compute the matrix of the bilinear form (\ref{pairing1}) in the
basis given by the columns of the matrix $\Phi_0$.

\begin{theorem} In the basis of the columns of the matrix $\Phi_0$ the bilinear
form (\ref{pairing1}) is given by the matrix
\beq\label{gram-0}
\left( \Phi_0^i(-\nu), \Phi_0^j(\nu)\right) 
=\left( q^{1/2} + q^{-1/2}\right)^2 \left[\eta\left(
 q^{1/2} e^{-\pi\, i\, R} e^{-\pi \, i\,
\hat\mu} +q^{-1/2} e^{\pi\, i\, R} e^{\pi \, i\,
\hat\mu}\right) \right]^{-1}_{ij}.
\eeq
\end{theorem}

\pf From the proof of the previous Theorem we obtained that
$$
\Phi_0(\nu) =\Phi(\nu) \, M(\nu)
$$
where
$$
M(\nu) = (1+q) (q\, S +S^T)^{-1} C.
$$
According to the formula (\ref{q-pair}) it remains to compute the product of the
following matrices
$$
M^T(-\nu) \left( q^{1/2} S + q^{-1/2} S^T\right) M(\nu).
$$
This can be easily done using
\eqa\label{c-s}
&&
S=C  e^{-\pi\, i\, R} e^{-\pi \, i\,
\hat\mu} \eta^{-1} C^T
\\
&&
S^T=C  e^{\pi\, i\, R} e^{\pi \, i\,
\hat\mu} \eta^{-1} C^T.
\label{s-s-t}
\eeqa
The Theorem is proved.\epf

In particular, for $\nu=0$ (i.e., $q=1$) one obtains the Gram matrix of the
intersection form in the basis given by the flat coordinates with good behaviour
at $\lambda=\infty$. Recall that this Gram matrix was used in the free field
realization of the Virasoro operators associated with the given Frobenius
manifold \cite{normal}. The general formula (\ref{gram-0}) was used in
\cite{normal} for regularization of the Virasoro operators in the resonant case.

\setcounter{equation}{0}
\setcounter{theorem}{0}

\section{Examples and applications}\label{sec4}\par

\subsection{Almost duality in the singularity theory}\par

Let us compute the twisted periods  for the Frobenius
structure arising, according to the K.Saito theory \cite{saito} of primitive forms (see
also \cite{prim,msaito,oda,manin,sabbah,hertling}), on the base of the universal unfolding
(i.e., on the parameter space of a versal deformation of the singularity)
 of an isolated
singularity $f(x)$,    
$x\in B\subset {\mathbb C}^N$ for a sufficiently small ball $B$, $f(x)=0$, 
$df(0)=0$.  Denote $f_t(x)$ the corresponding versal
deformation, $t=(t^1, \dots, t^n)$ are the {\it flat coordinates} \cite{prim,
hertling} on the
base $M_f$ of the versal deformation, $n$ is the Milnor number of the 
singularity. The discriminant $\Sigma\subset M_f$ consists of those values
of the parameters $t$ for which the fiber $f_t^{-1}(0)\cap B$ is singular.
The locus $\Sigma_\lambda$ (see (\ref{discr}) is defined in a similar way by the
conditions of singularity of the fiber $f_t^{-1}(\lambda)\cap B$ for $\lambda$ 
sufficiently close to $0$.
The {\it period mapping} in the singularity theory \cite{prim, oda}
\eqa\label{per-map}
&&
\pi : M_f \setminus \Sigma   \to H^{N-1}\left( f_t^{-1}(0)\cap B\right)
\nn\\
&&
t \mapsto [ \omega_t(x) ]
\eeqa
is obtained by choosing a holomorphic $(N-1)$-form on
$\left(M_f\setminus\Sigma\right) \times B$ closed along the fibers
$f_t^{-1}(0)$. The coordinate representation of the period mapping is a
multivalued vector-function
\beq\label{periods}
\pi(t)= \left( \oint_{\sigma_1}\omega_t(x), \dots,
\oint_{\sigma_n}\omega_t(x)\right)
\eeq
where $\sigma_1$, \dots, 
$\sigma_n\in \tilde H_{N-1}\left( f_t^{-1}(0)\cap B; {\mathbb
Z}\right)$ is a basis of vanishing cycles. Multivaluedness of the period mapping
is encoded by an action 
of the monodromy group $W$ of the singularity in the space of vanishing 
homologies. If the differential form
$\omega_t(x)$ is chosen in a clever way then $\pi$ is a local diffeomorphism.
Then the isomorphism
\beq\label{dif-per}
d\pi_*: H_{N-1}\left( f_t^{-1}(0)\cap B; {\mathbb
Z}\right)\to T_t^*M_f
\eeq
dual to the differential $d\pi$ defines a bilinear form $(~,~)$ on the cotangent
bundle to $M_f\setminus \Sigma$ induced by the intersection index pairing on the
homologies. This bilinear form is symmetric/skew-symmetric for $N-1$ even/odd.

To identify the above period mapping and intersection form of the singularity
theory with those coming from the theory of Frobenius manifolds we are to assume
that: 1). $N\equiv 1 ({\rm mod}\, 4)$. 2). The differential form $\omega_t(x)$ 
must
be a good primitive form in the sense of the K.Saito's theory of primitive forms
\cite{saito}. For the case of simple singularities an explicit construction
of a good primitive form was obtained by M.Noumi \cite{noumi}. For a general hypersurface
singularity existence of a good primitive form has been proved by M.Saito
\cite{msaito}.
Recall \cite{lo} that for the case of simple singularities the period mapping ($N-1$
must be even) produces
an analytic isomorphism
\beq\label{per-simple}
M_f\to {\mathbb C}^n/W.
\eeq
The flat coordinates on $M_f$ in this case are given by a certain remarkable
basis of homogeneous $W$-invariant polynomials. Their intrinsic contruction in
terms of the Weyl group and its generalization to an arbitrary finite Coxeter
group was found in \cite{sys}. The construction of the Frobenius structure on the orbit
spaces of finite Coxeter groups was obtained in \cite{cox} (see also
\cite{cime}). Some further
generalizations of this construction see in \cite{dz}.

We will consider here only the case of simple singularities labelled 
by Weyl groups of simply-laced Lie algebras \cite{arn}, i.e., by
Dynkin diagrams of $ADE$ type.
The reader may have in mind the example of the $A_n$ singularity
\beq\label{a-n}
f(x) =x^{n+1}, ~~f_t(x) =x^{n+1} + a_1 x^{n-1}+\dots +a_n , ~~x\in {\mathbb C}.
\eeq
Here $a_1$, \dots, $a_n$ are polynomials of $t^1, \dots, t^n$ that are
constructed in the following way \cite{sys}. Let us consider the series
$$
k:= f_t^{1\over n+1}(x) = x + {a_1\over n+1} x^{-1} +O\left( x^{-2}\right).
$$
The flat coordinates $t^\alpha=t^\alpha(a_1, \dots, a_n)$ are defined as the
first $n$ nontrivial coefficients of the inverse function expansion
$$
x = k +{1\over n+1} \left( {t^n\over k} + {t^{n-1}\over k^2} + \dots + {t^1\over
k^n}\right) + O\left( k^{-(n+2)}\right).
$$
The discriminant $\Sigma$ consists of all polynomials
with multiple roots. The subspace $t\in M_{s\, s}\subset M_f$ consists of all
polynomials $f_t(x)$ that are good Morse functions, i.e., all their critical
points are nondegenerate and the critical
values are pairwise distinct.
The dependence $f_t(x)$ on the 
flat coordinates
satisfies the following identities \cite{eguchi}
\eqa\label{eguchi}
&&
\phi_\alpha\phi_\beta = c_{\alpha\beta}^\gamma \phi_\gamma +
K_{\alpha\beta}^a {\partial f_t\over \partial x^a}
\nn\\
&&
\partial_\alpha\phi_\beta ={\partial K_{\alpha\beta}^a\over \partial x^a}.
\eeqa
Here 
\beq\label{phi}
\phi_\alpha=\phi_\alpha(x;t):=\partial_\alpha f_t,
\eeq
$c_{\alpha\beta}^\gamma=c_{\alpha\beta}^\gamma(t)$ are the structure constants
tensor of the Frobenius manifold,
the polynomials $K_{\alpha\beta}^a=K_{\alpha\beta}^a(x;t)$ are defined by
(\ref{eguchi}),
we assume also a summation w.r.t. Latin indices $a$ etc. from $1$ up to $N$.
The versal deformation can be chosen to be a quasihomogeneous one, i.e., it
satisfies also the identity
\beq\label{q-h-v}
f_t  = \sum_{a=1}^N r_a x^a {\partial f_t\over \partial x^a} + \sum_\alpha
(1-q_\alpha) t^\alpha \phi_\alpha
\eeq
with some rational numbers $r_1$, \dots, $r_N$ satisfying
\beq\label{q-h-r}
\sum_ar_a={N-d\over 2},
\eeq
\beq\label{q=h-d}
d=1-{2\over h}, ~~q_\alpha=1-{m_\alpha+1\over h},
\eeq
$h$ being the Coxeter number and $m_1$, \dots, $m_n$ the exponents of the Dynkin
diagram of the singularity.

We want to show that  the twisted periods can be
computed as the loop integrals of the following form
\beq\label{twist}
\tilde p(t; \nu) =\oint_\gamma f_t^{\nu-{N-1\over 2}}(x)d^Nx.
\eeq
Here $\gamma$ is a cycle in $H_N\left( B\setminus f_t^{-1}(0), {\bf
L}(q)\right)$ where the local system ${\bf L}(q)$ is defined by multiplication
by $(-1)^{N-1\over 2}q$ where $q:= e^{2\pi\, i\, \nu}$ (see details in
\cite{giv}). To this end we
are to prove that 
the (multivalued) functions $\xi_\alpha(t;\nu)$ of $t^1$, \dots, $t^n$ 
\beq\label{twist1}
\xi_\alpha=\oint_\gamma \phi_\alpha(x;t) f_t^{\nu'-1}(x)\, d^Nx,
\eeq
$$
\nu':=\nu-{N-1\over 2}
$$
satisfy the system (\ref{g-m-dual})
We will omit the reference to the cycle $\gamma$ in the computations. What we
will use of the symbol of loop integral is the usual properties  that the integral of a total
derivative vanishes and that the derivatives of the integrals along the parameters $t$ coincide
with the integrals of the derivatives.

\begin{prop}  The functions (\ref{twist1}) satisfy the system (\ref{g-m-dual}) for the
gradients of the twisted periods of the Frobenius manifold $M_f$.
\end{prop}

\pf Multiplying (\ref{q-h-v}) by $\phi_\beta$ and using (\ref{eguchi}) we obtain 
$$
\sum _a r_a x^a \phi_\beta {\partial f_t\over \partial x^a}+
{\mathcal U}_\beta^\gamma
\phi_\gamma +\sum_\alpha(1-q_\alpha) t^\alpha K_{\alpha\beta}^a {\partial
f_t\over \partial x^a}= \phi_\beta f_t(x).
$$
Another identity we obtain differentiating (\ref{q-h-v}) along $\partial_\beta$:
\beq\label{quasi-phi}
\sum_a r_a x^a {\partial \phi_\beta\over \partial x^a} +\sum_\alpha (1-q_\alpha)
t^\alpha {\partial K_{\alpha\beta} ^a \over \partial x^a} = q_\beta \phi_\beta.
\eeq
Now, differentiating (\ref{twist1}) and using (\ref{eguchi}) we obtain
$$
\partial_\epsilon\xi_\beta =(\nu'-1) c_{\epsilon\beta}^\gamma \oint \phi_\gamma
f_t^{\nu'-2}d^Nx.
$$
Here we have eliminated the divergence term
$$
{\pal K^a_{\epsilon\beta}\over \pal x^a} f_t^{\nu'-1} + (\nu'-1)
K_{\epsilon\beta}^a {\pal f_t\over \pal x^a } f_t^{\nu'-2}={\pal\over \pal x^a}
\left[ K^a_{\epsilon\beta} f_t^{\nu'-1}\right].
$$
Multiplying the last equation  by ${\mathcal U}_\alpha^\epsilon$ and 
using associativity and (\ref{eguchi}), (\ref{q-h-v}) we obtain
$$
{\mathcal U}_\alpha^\epsilon \partial_\beta\xi_\epsilon =(\nu'-1)
c_{\alpha\beta}^\gamma \xi_\gamma - c_{\alpha\beta}^\gamma
\oint \sum_a\left[ \phi_\gamma  r_a x^a 
+ \sum_\epsilon (1-q_\epsilon) t^\epsilon K_{\epsilon\gamma}^a
\right] {\partial 
f_t^{\nu'-1}\over \partial x^a}
\, d^Nx.
$$
Integrating the last integral by parts and using (\ref{quasi-phi}) we arrive at
$$
{\mathcal U}_\alpha^\epsilon \partial_\beta\xi_\epsilon =\sum_\gamma(\nu'-1+\sum_ar_a+q_\gamma)
c_{\alpha\beta}^\gamma \xi_\gamma.
$$
Due to (\ref{quasi-phi}), this coincides with (\ref{g-m-dual}). Proposition is proved.

One can
construct a basis of vanishing cycles $\sigma_1$, \dots, $\sigma_n$ in
$\tilde H_{N-1}(f_t^{-1}(0)\cap B, {\mathbb Z})$ for any 
$t\in M_f\setminus \Sigma$. The basis varies continuously with small variations
of $t$.
It is of particular convenience to choose a so-called {\it distinguished basis}
of vanishing cycles corresponding to the properly ordered critical values
$u_1(t)$, \dots, $u_n(t)$ of $f_t(x)$ connected by non-intersecting paths to the
origin (see, e.g., \cite{agv}).
The Givental's
construction \cite{giv} gives a way to associate to it a basis $\gamma_1$, \dots,
$\gamma_n$ in  the homology $H_N\left( B\setminus f_t^{-1}(0), {\bf
L}(q)\right)$. Taking the integrals
\beq\label{twist2}
\tilde p_1(t;\nu)=\oint_{\gamma_1} f_t^{\nu'}(x)d^Nx, \dots, 
\tilde p_n(t;\nu)=
\oint_{\gamma_n} f_t^{\nu'}(x) d^Nx
\eeq
we obtain, locally on $M_f\setminus \Sigma$, the twisted
period mapping. Globally the monodromy of the mapping is described by
{\it twisted Picard - Lefschetz formulae} found in \cite{giv}). In the next
Section we will derive  the analogue of these formulae for the
monodromy of the twisted period mapping on an arbitrary semisimple Frobenius manifold.

\subsection{Frobenius and almost Frobenius structures on orbit spaces of
finite Coxeter groups}\par 

Let $W\subset End({\mathbb R}^n)$ be a finite irreducible
Coxeter group. In \cite{cox} it was constructed a structure of polynomial
Frobenius manifold on the orbit space
\beq\label{orb}
M={\mathbb C}^n / W.
\eeq
A coordinate system on the orbit space is given by choosing $n$
homogeneous $W$-invariant polynomials $y^1(x)$, \dots, $y^n(x)$
generating the ring ${\mathbb C}[x^1, \dots, x^n]^W$ of $W$-invariant
polynomials on ${\mathbb C}^n$. The first metric $(~,~)$ on the orbit space
reads
\beq\label{arnold}
\left( dy^i, dy^j\right) =\sum_{a,\,b}{\partial y^i\over \partial x^a}
{\partial y^j\over \partial x^b} G^{ab} = g^{ij}(y).
\eeq
Here $G^{ab}=\left( dx^a, dx^b\right)$ are the contravariant components
of a $W$-invariant Euclidean metric on ${\mathbb R}^n$. The components
$g^{ij}(y)$ are polynomials in $y^1$, \dots, $y^n$ (cf. \cite{sys}). The second
metric \cite{sys, saito} is given by
\beq\label{saito}
<dy^i, dy^j> ={\partial g^{ij}(y)\over \partial y^1}
\eeq
assuming that $h={\rm deg}\, y^1(x)$ is the maximum of the degrees of the
basic invariant polynomials. The discriminant $\Sigma\subset M$ consists
of all orbits containing less than $|W|$ points. Alternatively it can be
described as the image of the (complexified) mirrors of all reflections in
${\mathbb R}^n \subset {\mathbb C}^n$ w.r.t. the natural projection
$$
{\mathbb C}^n\to {\mathbb C}^n /W = M.
$$
Outside these mirrors the projection is a local diffeomorphism onto
$M\setminus\Sigma$. The period mapping 
$$
p^1(y^1, \dots, y^n), \dots, p^n(y^1, \dots, y^n)
$$
is given by inverting this local
diffeomorphism, i.e., by solving the system of algebraic equations
$$
y^1(p^1, \dots, p^n)=y^1,~ y^2(p^1, \dots, p^n)=y^2, \dots, y^n(p^1, \dots,
p^n)=y^n.
$$
The flat coordinates of the flat pencil (\ref{arnold}), (\ref{saito}) are determined by the system
\beq\label{f-c-lambda}
y^1(p^1, \dots, p^n)=y^1-\lambda,~ y^2(p^1, \dots, p^n)=y^2, \dots,
y^n(p^1, \dots,
p^n)=y^n
\eeq
where, we recall, the degree of the polynomial $y^1(x)$ is the maximal
one.

We give now the formula for the dual potential $F_*(x)$ for the case of
polynomial Frobenius manifolds defined on the orbit spaces of finite
Coxeter groups (see above). Let $\Delta_+$ be the set of all positive
roots of the Coxeter group $W\subset End({\mathbb R}^n)$. The hyperplanes
$$
(\alpha,x)=0, ~~\alpha\in \Delta_+
$$
are all the mirrors of the reflections in $W$. Let us normalize the roots
by the condition
$$
(\alpha,\alpha)=2.
$$
We will identify the roots $\alpha$ with the linear functions $x\mapsto
(\alpha,x)$.

\begin{theorem}  For any finite irreducible Coxeter group the function
$F_*(x)$ defined on the universal covering of 
\beq\label{m-star}
{\mathbb C}^n\setminus \cup_{\alpha\in\Delta_+} \left\{ (\alpha,x)=0\right\}
\eeq
has the form
\beq\label{star-cox}
F_*(x)={h\over 4} \sum _{\alpha\in \Delta_+} \alpha^2 \log \alpha^2.
\eeq
\end{theorem}

\pf For the Frobenius manifold under consideration the discriminant $\Sigma$
is a {\it finite} union of hyperplane in the Euclidean space ${\mathbb R}^n$.
From the Theorem \ref{c-th} we know that all the singularities of $F_*(x)$ must
be on these hyperplanes only. The third derivatives of $F_*(x)$ must have
singularities of the from (\ref{c-beh}). So, they are meromorphic functions
on ${\mathbb R}^n$ with simple poles along the mirrors of the reflections of the
Coxeter group.
Clearly the function $F_*(x)$ is
determined by these analytic properties uniquely up to adding of a quadratic
factor.
Let us check that the formula (\ref{m-star}) satisfies the needed analytic
properties.

Let $e_1$, \dots, $e_n$ be a basis of ${\mathbb R}^n$. As above we
denote
$$
G_{ab} =(e_a, e_b), ~~ \left( G^{ab}\right) = \left( G_{ab}\right)^{-1}.
$$
We also put
$$
\al_a:= (\alpha, e_a), ~~ a=1, \dots, n, ~~\al\in \Delta_+.
$$
Taking the triple derivatives of $F_*(x)$ we obtain
\beq\label{c-star-cox}
\cstar_{ab}^c(x) ={h\over 2} \sum_{\al\in\Delta_+} 
{\al_a \al_b \al_d\over (\al,
x)} G^{dc}.
\eeq
The singular part of this formula near the mirror $(\al,x)=0$ coincides with
(\ref{c-beh})
since
$$
{1\over 1-d} = {h\over 2}.
$$ 
The Theorem is proved.\epf

{\bf Remark 2.} Associativity of the family of algebras follows from Corollary \ref{cor1-2}. It can
be proved also by straightforward computation \cite{mg, ves0}. Remarkably, A.Veselov
recently found \cite{ves} other examples of solutions to the associativity
equations given by
a formula similar to (\ref{c-star-cox}). In the Veselov's examples the so-called
deformed root systems are used. Recall that deformed root systems were
discovered in the theory of multidimensional integrable linear differential
operators (see \cite{cfv} and references therein).We do not know if Veselov's
structures satisfy the Axiom {\bf AFM3}.

\begin{cor} The twisted periods as functions of the Euclidean
coordinates $x^1$, \dots, $x^n$ satisfy the following system of linear
differential equations with rational coefficients
\beq\label{twist-nu}
\pal_a\xi_b = {h\,\nu\over 2} \sum_{\al\in\Delta_+} 
{\al_a \al_b \al_c\over (\al,
x)} G^{cd}\xi_d, ~~ \xi_a = {\pal\over\pal x^a} \tilde p (t(x); \nu).
\eeq
\end{cor}

{\bf Example 2.} To construct the almost Frobenius structure for 
$W=W(A_n)$ it is convenient to start with the standard action of $W=S_{n+1}$
on ${\mathbb R}^{n+1}$ by permutations of the coordinates $x_0$, $x_1$, \dots,
$x_n$. Introduce the function
\beq\label{prepo}
F_*(x) ={n+1\over 8} \sum_{i<j} (x_i-x_j)^2 \log (x_i-x_j)^2.
\eeq
Together with the Euclidean metric
\beq\label{evklid}
(dx_i, dx_j) = \delta_{ij}
\eeq
the third derivatives of $F_*(x)$ give the following multiplication law
of tangent vectors $\pal_i := \pal / \pal x_i$
\beq\label{law}
\pal_i * \pal_j =-{n+1\over 2} {\pal_i - \pal_j\over x_i -x_j}, \quad i\neq j
\eeq
$$
\pal_i * \pal_i = - \sum_{j\neq i} \pal_i * \pal_j.
$$
The vector field 
$$
\sum_i \pal_i
$$
has zero products with everything. Factorizing over this direction one obtains
an almost Frobenius structure on
$$
M_* = \{ x_0 + x_1 + \dots + x_n =0 \} \setminus \cup _{i<j} \{ x_i =x_j\} .
$$
The unity = Euler vector field equals
$$
E={1\over n+1} \sum_i x_i \pal_i
$$
the vector field $e$ of the axiom {\bf AFM3} reads
$$
e=-\sum_i {1\over f'(x_i)} \pal_i, \qquad f(x):= \prod_{i=0}^n (x-x_i).
$$
The equation (\ref{def-star}) for the twisted periods coincides
with the classical Euler - Poisson - Darboux equations
\beq\label{eul-darb}
\pal_i \pal_j \tilde p = - {\nu\over x_i - x_j} (\pal_i \tilde p - \pal_j \tilde
p), \quad i\neq j.
\eeq 
Recently these equations proved to play an important role in the theory of
dispersive waves \cite{tian, grava}.

\subsection{Almost duality and Shephard groups}\par

Shephard groups are the symmetry groups of regular complex polytopes. It is
a subclass of finite {\it unitary reflection groups}. By definition
a unitary reflection is a linear transformation
$$
g: {\mathbb C}^n\to {\mathbb C}^n
$$
having a hyperplane of fixed points such that the only nontrivial eigenvalue of
$g$ is a root of unity. A finite unitary reflection group $G$ by definition is a finite
subgroup in $GL(n,{\mathbb C})$ generated by unitary reflections. The book
\cite{orlik} is an excellent introduction into the theory of unitary reflection
groups.

Of course, any finite Coxeter group is a Shephard group.  
Besides this, there
are two infinite series and fifteen exceptional cases of irreducible 
Shephard groups (see Table 1 below). It was discovered in \cite{os1} that, for any
Shephard group $G$ there is an accompanying finite Coxeter group $W$. The group
$G$ is uniquely determined by the pair $(W,\kappa)$ where the number $\kappa$
was defined in (\ref{order}). Let us represent this correspondence
describing the Shephard group associated with $(W,\kappa)$ in terms of
the monodromy of twisted periods on the Frobenius manifold $M_W$.

According to Chevalley theorem, the quotient
\beq\label{mg}
M_G ={\mathbb C}^n /G
\eeq
carries a natural structure of graded affine variety. One can use a basis
homogeneous
$G$-invariant polynomials $f^1(z)$, \dots, $f^n(z)$, $z=(z^1, \dots, z^n)\in
{\mathbb C}^n$ as the graded affine coordinates on $M_G$. Let us order them such
that
\beq\label{order}
\deg f^1(z) = {\rm max} = h_G > \deg f^2(z) \geq \dots > \deg f^n(z) ={\rm
min}=:2 \kappa.
\eeq
The class of Shephard groups is selected by the following remarkable property
\cite{os1}, see also \cite{orlik}.

\begin{theorem} For a Shephard group the Hessian of the basic
invariant polynomial of the lowest degree is a nondegenerate matrix 
for generic $z$.
Conversely,
nondegeneracy of the Hessian completely
characterizes irreducible Shephard groups among all unitary reflection groups.
\end{theorem}

We give here a nice differential-geometric interpretation of the Hessian
reinterpreting the results of Orlik and Solomon \cite{os1, os2}.

Let us denote
\beq\label{gij0}
h_{ij}(z) = {\pal^2 f^n(z)\over \pal z^i \pal z^j}, \quad \left(
h^{ij}(z)\right) = \left( h_{ij}(z)\right)^{-1}.
\eeq

The $G$-bilinear form 
\beq\label{gij1}
(df^i, df^j):= {\pal f^i\over \pal z^k} {\pal f^j\over \pal z^k} h^{kl}(z)
\eeq
is defined on the cotangent planes $T^* M_G$ everywhere due to the following
statement \cite{os1}.

\begin{theorem} The functions $(df^i, df^j)$ are polynomials. The determinant
$\det (df^i, df^j)$ vanishes exactly on the collection of mirrors of the group
$G$.
\end{theorem}

Because of $G$-invariance one can represent $(df^i, df^j)$ as polynomials in
$f^1(z)$, \dots, $f^n(z)$
\beq\label{gij2}
(df^i, df^j)=:g^{ij}(f), \quad f=(f^1, \dots, f^n).
\eeq

We obtain a polynomial flat metric on $T^*M_G$. Define an analogue of K.Saito
metric by
\beq\label{gij3}
\eta^{ij}(f):= {\pal g^{ij}(f)\over \pal f^1}.
\eeq

From \cite{os2} it follows

\begin{theorem} The metric $(\ref{gij3})$ on
$T^*M_G$ does not degenerate anywhere. 
\end{theorem}

The main step in establishing a connection between Shephard groups and Coxeter
groups is

\begin{theorem} \label{main-sh} The polynomial metrics (\ref{gij2}), (\ref{gij3})
together with the unity vector field
$$
e={\pal\over\pal f^1}
$$
and the Euler vector field
$$
E= h_G^{-1} \sum_{i=1}^n (\deg f^i) f^i {\pal \over \pal f^i}
$$
form a flat pencil.
\end{theorem}

The flat coordinates for the metric (\ref{gij3}) are given by a distinguished
system of homogeneous {\it flat generators} in the ring of $G$-invariant
polynomials. Flat generators exist also for other unitary reflection groups
\cite{orlik}, but
they are not flat coordinates of a natural metric on the orbit space if $G$ is
not a Shephard group.

The following two statements immediately follow from the Theorem \ref{main-sh}.

\begin{cor} The orbit space $M_G$ of a Shephard group $G$ carries a natural
polynomial Frobenius structure.
\end{cor}

Due to the Hertling's theorem \cite{hertling} the Frobenius manifold $M_G$
must be isomorphic to the orbit space of a finite irreducible Coxeter group
$W$. This is just the accompanying Coxeter group for $G$!

\begin{cor} The
Hessian quadratic form
\beq\label{hess}
ds^2:= {\pal^2 f^n(z)\over \pal z^i \pal z^j} dz^i dz^j
\eeq
is a flat metric on
\beq\label{g-mir}
M_G^* := {\mathbb C}^n \setminus \{ {\rm mirrors} ~ {\rm of}~ G \}.
\eeq
\end{cor}

If $G$ is itself a Coxeter group then the metric (\ref{hess})
coincides, up to a constant factor, with the Killing form.

Let us make a digression about the flat metrics representable in the
Hessian form. Remarkably, our old friend associativity equations arises
also in this setting!

\begin{prop} \cite{kito}. Let $f(z)$ be a smooth function of $z=(z^1, \dots, z^n)$
such that the Hessian does not degenerate in some domain $\Omega\subset {\mathbb
R}^n$
$$
\det \left( {\pal ^2 f(z)\over \pal z^i \pal z^j}\right) \neq 0.
$$
Denote
\beq\label{h-inv}
\left( h^{ij}(z)\right) := \left( {\pal ^2 f(z)\over \pal z^i \pal
z^j}\right)^{-1}
\eeq
the inverse matrix. Then the metric
$$
ds^2:= {\pal^2 f(z)\over \pal z^i \pal z^j} dz^i dz^j
$$
on $\Omega$ has a zero curvature {\rm iff} the function $f$ satisfies the
following system of associativity equations
\beq\label{h-ass}
{\pal^3 f(z)\over \pal z^i \pal z^j \pal z^s} h^{st}(z) 
{\pal^3 f(z)\over \pal z^t \pal z^k \pal z^l}=
{\pal^3 f(z)\over \pal z^l \pal z^j \pal z^s} h^{st}(z) 
{\pal^3 f(z)\over \pal z^t \pal z^k \pal z^i}
\eeq
for all $i$, $j$, $k$, $l$.
\end{prop}

\pf Let us denote by superscripts the partial derivatives of $f$ w.r.t. $z^i$,
$z^j$, etc. An easy calculation gives the Christoffel coefficients for the
Levi-Civita connection for the metric $ds^2$:
$$
\Gamma_{ij}^k = {1\over 2} h^{ks} f_{sij}.
$$
From this one readily derives the following formula for the Riemann curvature
tensor of the metric
$$
R_{ijl}^k = {1\over 4} h^{kp} f_{pqj} h^{qs} f_{sil} -{1\over 4} h^{kp}
f_{pqi}h^{qs} f_{slj}.
$$
Vanishing of this expression is equivalent to (\ref{h-ass}).\epf

We will now describe a natural class of flat Hessian metrics associated with
a Frobenius manifold.

\begin{prop} Let $M$ be an arbitrary Frobenius manifold with the charge $d\neq
1$. Let
\beq\label{z-coor}
z^i := \tilde p^i(t; \nu), ~~i=1, \dots, n
\eeq
be a system of independent twisted periods on $M\setminus\Sigma$. 
Here $\nu$ is a fixed complex number satisfying
\beq\label{no-nu}
\nu\neq {1-d\over 2}.
\eeq
Denote
\beq\label{z-met}
g_{ij}(z):= \left({\pal\over \pal z^i}, {\pal\over \pal z^j}\right)
\eeq
the Gram matrix of the intersection form (\ref{int-form}) on $TM$ written
in the coordinates $z^1$, \dots, $z^n$. Put
\beq\label{f-z}
f(z)=\left[ {1-d\over 2} -\nu\right]\,t_1
\eeq
where $t_1=\eta_{1\al}t^\al$ must be represented as a function of the
coordinates $z$.
Then this flat metric can be written
in the Hessian form
\beq\label{hess-form}
g_{ij}(z) = {\pal ^2 f(z)\over \pal z^i \pal z^j}.
\eeq
\end{prop}

Observe that the function $f(z)$ is a homogeneous function of $z=(z^1, \dots,
z^n)$ of the degree
\beq\label{deg-f}
\deg f(z) ={1-d\over {1-d\over 2} + \nu}.
\eeq
This easily follows from 
$$
Lie_E t_1 = (1-d) t_1
$$
and
$$
Lie_E z^i = {1-d\over 2} + \nu
$$
(see (\ref{quasi-p})). 

\pf Let us compute the Hessian of the function (\ref{f-z}). We have, due to 
(\ref{t1})
\beq\label{eq1}
{2\over 1-d} {\pal^2 f\over \pal z^i \pal z^j} 
= G_{ab} {\pal p^a\over \pal z^i}{\pal p^b\over \pal z^j}
+ G_{ab} p^a {\pal^2 p^b\over \pal z^i \pal z^j}.
\eeq
The first term coincides with the metric (\ref{z-met}). Let us compute the
second one.

Rewriting
$$
{\pal^2 p^b\over \pal z^i \pal z^j}={\pal\over \pal p^c} \left( {\pal p^b\over
\pal z^i} \right) {\pal p^c\over \pal z^j}
$$
and using
$$
{\pal p^b\over
\pal z^s } {\pal z^s\over \pal p^d} = \delta^b_d
$$
yields
$$
{\pal^2 p^b\over \pal z^i \pal z^j}=-{\pal p^b\over \pal z^s} {\pal^2 z^s\over
\pal p^c \pal p^d} {\pal p^d\over \pal z^i} {\pal p^c\over \pal z^j}.
$$
We now use the equation (\ref{def-star}) for twisted periods to recast the
second term in (\ref{eq1}) as follows
$$
G_{ab} p^a {\pal^2 p^b\over \pal z^i \pal z^j}= -\nu\, p^a \cstar_{acd} {\pal
p^c\over \pal z^i} {\pal p^d\over \pal z^j}.
$$
The last step in the proof is to use is that, the vector field
$$
E={1-d\over 2} p^a {\pal\over\pal p^a}
$$
is the unity on the dual almost Frobenius manifold $M_*$. Therefore
$$
p^a \cstar_{acd} ={2\over 1-d}G_{cd}.
$$
This proves that the second term in (\ref{eq1}) is proportional to the first
one:
$$
G_{ab} p^a {\pal^2 p^b\over \pal z^i \pal z^j}= -{2\nu\over 1-d}G_{cd}{\pal
p^c\over \pal z^i} {\pal p^d\over \pal z^j}.
$$
This proves the Proposition.\epf

Summarizing the above results we arrive at the main statement of this Section.

\begin{theorem} Let $W$ be a finite Coxeter group acting in the $n$-dimensional
space. Denote $h$ the Coxeter number of $W$.
Let $e_1$, \dots, $e_n$ be a basis of simple roots normalized by $(e_i, e_i)=2$.
Introduce an upper triangular matrix $S$ such that
$$
S_{ii}=1, \quad S_{ij}=(e_i, e_j) \quad {\rm for} ~ i<j.
$$
Let $\nu$ be a rational number such that
$$
\det (q S+S^T) \neq 0, \quad q:= e^{2\pi\, i}
$$
and the monodromy matrices
$M_1(q)$, \dots, $M_n(q)$ of the form (\ref{mono-stokes}) generate a finite
irreducible subgroup $G$ in $GL(n, {\mathbb C})$. Then $G$ is a Shephard group
with the accompanying Coxeter group $W$ and
\beq\label{kap}
\kappa = {1\over 1-h \, \nu}
\eeq
where $h$ is the Coxeter number of $W$. Conversely, if $(W,\kappa)$ is as in
Table 1 then the monodromy matrices $M_1(q)$, \dots, $M_n(q)$
of twisted periods $\tilde p^a(t; \nu)$ with
\beq\label{nu}
\nu = {1\over h} \left( 1-{1\over \kappa}\right)
\eeq
generate the Shephard group associated with $(W,\kappa)$.
\end{theorem}

\def\qq{&\qquad}
$$
\begin{matrix}
W\qq \ \kappa \qq \  G \\ \\ \\
A_1 \qq \ {r\over 2} \qq \ C(r) \\ \\
B_n \qq \ {r\over 2} \qq \ G(r,1,n) \\ \\
A_2 \qq \ 2 \qq \ G_4 \\
A_2 \qq \ 4 \qq \ G_8 \\
A_2 \qq \  10 \qq \ G_{16} \\
B_2 \qq \ 3  \qq \  G_5 \\
B_2 \qq \ 6 \qq \ G_{10} \\
B_2 \qq \ 15 \qq \ G_{18} \\
I_2(5) \qq \ 6 \qq \ G_{20}\\
G_2 \qq \ 2 \qq \  G_6\\
G_2 \qq \ 4 \qq \ G_9 \\
G_2 \qq \ 10\qq \ G_{17}\\
I_2(8) \qq \ 3\qq \ G_{14}\\
I_2(10) \qq \ 6\qq \ G_{21}\\
A_3 \qq \ 3\qq \ G_{25}\\
B_3 \qq \ 3\qq \ G_{26}\\
A_4 \qq \ 6\qq \ G_{32}
\end{matrix}
$$
\vskip 0.5cm
$$
\text {Table 1}
$$
\vskip 0.5cm
\vskip 0.5cm
 
{\bf Example 3.} Let us consider the particular case of the group $G_{25}\subset
GL(3)$. This
is the group of symmetries of the celebrated Hessian configuration consisting of
9 inflections of a generic plane cubic. The ring of the $G_{25}$-invariant
polynomials is generated by the classical Maschke polynomials
\eqa\label{maschke}
&&
C_6(z) = z_1^6 + z_2^6 + z_3^6 - 10 (z_1^3 z_2^3 + z_2^3 z_3^3 + z_3^3 z_1^3)
\nn\\ 
&&
C_9(z) = (z_1^3 - z_2^3)(z_2^3-z_3^3)(z_3^3-z_1^3)
\nn\\
&&
C_{12}(z) = (z_1^3+z_2^3 + z_3^3) \left[ (z_1^3+z_2^3 + z_3^3)^3 + 216 z_1^3
z_2^3 z_3^3\right].
\eeqa
The Hessian flat metric reads
\eqa\label{he-ma}
&&
{1\over 30} ds^2 =3(z_1^4 dz_1^2 + z_2^4 dz_2^2 + z_3^4 dz_3^2)
-2(z_1^3+z_2^3+z_3^3) (z_1 dz_1^2 + z_2 dz_2^2 + z_3 dz_3^2) 
\nn\\
&&
\quad\quad\quad -6 (z_1^2 z_2^2 dz_1 dz_2+ z_2^2 z_3^2 dz_2 dz_3+z_3^2 z_1^2 dz_3
dz_1).
\eeqa
The flat coordinates $x_1$, $x_2$, $x_3$ of this metric are algebraic functions
of $z_1$, $z_2$, $z_3$ to be determined from the system
\eqa\label{alg-sys}
&&
C_6(z) = x_1^2 + x_2^2 + x_3^2
\nn\\
&&
C_9(z) = {1\over 4} x_1 x_2 x_3
\nn\\
&&
C_{12}(z) = (x_1^2 + x_2^2 + x_3^2)^2 - 3 (x_1^2 x_2^2 + x_1^2 x_3^2 + x_2^2
x_3^2),
\eeqa
$$
{\pal z^i\over \pal x^k} {\pal^2 C_6(z)\over \pal z^i \pal z^j} {\pal z^j\over
\pal x^l} = {10\over 3} \delta_{kl}.
$$
The inverse functions $z_i=z_i(x_1, x_2, x_3)$ are twisted periods
on the Frobenius manifold $M_{A_3}$ with $\nu= {1\over 6}$. Comparing with the
integral representation (\ref{twist}) of the twisted periods we obtain an
interesting realization of the group $G_{25}$ by the monodromy of Abelian
integrals of the form
$$
z=\oint (x^4 + a x^2 + bx + c)^{1\over 6} dx.
$$
The Theorem \ref{lef} gives the following matrix realization of the generators
of $G_{25}$
$$
M_1=\left( \begin{matrix} -q & q & 0 \\ 0 & 1 & 0 \\ 0 & 0 & 1
\end{matrix}\right), \quad M_2=\left( \begin{matrix} 1 & 0 & 0 \\ 1 & -q & q \\
0 & 0 & 1 \end{matrix}\right), \quad M_3=\left( \begin{matrix} 1 & 0 & 0 \\
0 & 1 & 0 \\ 0 & 1 & -q\end{matrix}\right)
$$
where
$$
q=e^{\pi \, i\over 3}.
$$

\subsection{Twisted periods for $QH^*(CP^4)$ and the mirror of the quintic}\par

Let $M=QH^*(CP^d)$ be the Frobenius manifold corresponding to the quantum
cohomology of the $d$-dimensional complex projective space. It is a semisimple
Frobenius manifold of the dimension $n=d+1$. The flat coordinates $t_1$, $t_2$,
\dots, $t_{d+1}$
on it are in
one-to-one correspondence with the standard basis $1$, $\omega$, \dots,
$\omega^d$ in $H^*(CP^d)$ (we now write all lower indices for the sake of
graphical simplicity). Here the two-form $\omega\in H^2(CP^d)$ is normalized by
the condition
$$
\int_{CP^d} \omega^d =1.
$$

The algebra structure on the tangent planes $T_tM$
becomes particularly simple at the points of the {\it small quantum cohomology
locus} 
\beq\label{small}
t=(t_1, t_2, 0, \dots, 0).
\eeq
At this points we have
\beq\label{quant}
T_tM\simeq {\mathbb C}[\omega]/\omega^{d+1} =e^{t_2}, ~~\omega\leftrightarrow
\partial_2.
\eeq
At the point of classical limit $t_2\to -\infty$ the algebra (\ref{quant}) goes to
the classical cohomology algebra of the projective space. Let us describe the
twisted periods at the locus (\ref{small}). It will be convenient to consider them
as a function of $t_2$ and $\lambda$ introducing $\lambda $ as in the beginning
of Section 3.

\begin{prop} The twisted periods $\tilde p =\tilde p(t; \nu)$ at the
points (\ref{small}) are determined from the following hypergeometric equation
\beq\label{cp-per}
\partial_2^{d+1} \tilde p =e^{t_2} t_1^{-(d+1)} \prod_{m=0}^d \left[ -(d+1)
\partial_2 +{1-d\over 2} +\nu -m\right] \, \tilde p
\eeq
and from the quasihomogeneity condition (\ref{quasi-p})
\beq\label{quasi-cp}
t_1 \partial_1 \tilde p + (d+1) \partial_2 \tilde p = \left( {1-d\over 2}
+\nu\right) \tilde p.
\eeq
\end{prop}

\pf The simplest way to derive (\ref{cp-per}) is to represent $\tilde p$ as the Laplace
integral 
$$
\tilde p = {i\over \sqrt{2\pi}} \int_0^\infty \tilde t(t; z) e^{-\lambda\, z}
{dz\over z^{\nu+{1\over 2}}}
$$
(cf. (\ref{lap0}))
and then to use the differential equation for the deformed flat
coordinates $\tilde t$ on the original Frobenius manifold. On the locus
(\ref{small})
 the latter
reads
\eqa\label{hyper-cp1}
&&
\partial_2^{d+1} \tilde t =z^{d+1} e^{t_2} \tilde t
\\
&&
\partial_1 \tilde t  = z \tilde t
\nn\\
&&
z\partial_z \tilde t  =t_1 \partial_1 \tilde t 
+ (d+1) \partial_2 \tilde t +{d-2\over 2} \tilde t.\label{hyper-cp2}
\eeqa
Substituting 
$$
\tilde t(t; z) = z^{\nu+{1\over 2}} \oint e^{\lambda\, z} \tilde p(t; \lambda)
$$
into the equation (\ref{hyper-cp2}) and integrating by parts we obtain
$$
\left\{ {1\over (d+1)^{d+1}} \left[ -\lambda \partial_\lambda +{1-d\over 2}
+\nu\right]^{d+1} -e^{t_2} \left( - \partial_\lambda\right)^{d+1}\right\} \,\tilde
p =0.
$$ 
We obtain a similar equation for the dependence of $\tilde p$ on $t_1$ since
$\partial_1 =-\partial_\lambda$. Setting then $\lambda$ to zero and using the
quasihomogeneity condition (\ref{quasi-cp})), we obtain (\ref{cp-per}). 
Proposition is proved.\epf

\begin{cor} Odd periods on the Frobenius
manifold $M=QH^*(CP^4)$ at the points (\ref{small}) with $t_1=-1$ are given by the periods of the holomorphic three-form
on the Calabi - Yau three-fold dual to the quintic in
$CP^4$.
\end{cor}

\pf The equation (\ref{cp-per}) for the
odd periods can be integrated once in $t_2$ to produce, at $t_1=-1$, the
Picard - Fuchs equation \cite{candelas} for the periods of the Calabi - Yau three-fold 
$$
u_0 + \dots + u_4 =1, ~~u_0 \dots u_4 =e^{t_2}
$$
dual
to the quintic in $CP^4$:
\beq\label{cand}
\partial_2^4 \tilde p = 5 \, (5\partial_2 +1) (5\partial_2+2) (5\partial_2+3) (5
\partial_2+4) \tilde p,
\eeq
$$
\tilde p(t_2) =\oint {du_0 \wedge \dots \wedge du_4 \over d(u_0+\dots
+u_4)\wedge d(u_0 \dots u_4)}.
$$
Corollary is proved.\epf

\subsection{From a Frobenius manifold to the Seiberg - Witten prepotential}\par

Let $M$ be an arbitrary $n$-dimensional Frobenius manifold with eigenvalues of
${\mathcal V}$ distinct from $1/2$. We will study the properties of the odd periods
on the $2n$-dimensional manifold
\beq\label{tensor}
M\otimes QH^*(CP^1).
\eeq

The construction of tensor product of Frobenius manifolds, introduced by
R.Kauf\-mann, M.Kontsevich, and Yu.I.Manin, generalizes the procedure of
computation of Gromov - Witten invariants of Cartesian product of two smooth
projective manifolds. Let us denote $t^1$, \dots, $t^n$ the flat coordinates on
$M$. For simplicity we will assume that ${\mathcal V}=\hat\mu$ is a diagonalizable matrix and we
choose the flat coordinates in such a way that
$$
\eta_{\alpha\beta} =\delta_{\alpha+\beta, n+1}, ~~{\mathcal V}={\rm diag}\, (\mu_1,
\dots, \mu_n), ~~\mu_\alpha+\mu_{n-\alpha+1}=0.
$$
The flat coordinates on $QH^*(CP^1)$ we redenote $(t^1, s)$, so the potential
and the Euler vector field read
\beq\label{cp1}
F_{CP^1}={1\over 2} (t^1)^2 s + e^s , ~~E=t^1\partial_1 + 2 \partial_s.
\eeq
The tangent spaces to (\ref{tensor}) has the structure of the tensor product
\beq\label{tensor1}
T\left( M\otimes QH^*(CP^1) \right)=TM \otimes H^* (CP^1).
\eeq
So the flat coordinates on the tensor product can be naturally labelled by pairs
$(\alpha', \alpha'')$ with $\alpha'=1, \dots, n$, $\alpha''=1, \, 2$. We
identify $t^{\alpha'\, 1''}$ with $t^\alpha$ and $t^{1'\, 2''}$ with $s$, and
we consider the ``coordinate cross''
\beq\label{cross}
t^{\alpha'\, 2''}=0, ~~\alpha'=2, \dots, n.
\eeq
The points of the coordinate cross will be coordinatized by pairs $(t, s)$,
$t=(t^1, \dots, t^n)\in M$. 

The Frobenius structure on (\ref{tensor}) is uniquely determined, according to 
\cite{kauf}, in a
neiborhood of the coordinate cross by the initial condition that, on the
coordinate cross, 
\beq\label{init}
T_{(t,s)}M\otimes QH^*(CP^1) =T_tM \otimes T_s QH^*(CP^1)
\eeq
is an isomorphism of Frobenius algebras. The Euler vector field of 
(\ref{tensor}) on the
coordinate cross reads
\beq\label{e-cross}
E_{M\otimes QH^*(CP^1)} =E+2\partial_s.
\eeq
The operator $\mu$ of (\ref{tensor}) thus acts as follows
\eqa\label{mu-tens}
&&
\mu(e_\alpha\otimes e_1) =(\mu_\alpha-{1\over 2}) e_\alpha\otimes e_1
\nn\\
&&
\mu(e_\alpha\otimes e_2) = (\mu_\alpha+{1\over 2}) e_\alpha\otimes e_2.
\eeqa
Hence the Poisson bracket (\ref{pb}) on the tensor product has the form
\eqa\label{pb-tens}
&&
\left\{ t^{\alpha' 1''}, t^{\beta'2''}\right\} =\eta^{\alpha\beta}
(\mu_\beta+{1\over 2})
\nn\\
&&
\left\{ t^{\alpha' 2''}, t^{\beta'1''}\right\} =\eta^{\alpha\beta}
(\mu_\beta-{1\over 2}),
\eeqa
other brackets vanish. This Poisson bracket induces a symplectic structure on the
tensor product (\ref{tensor})
\beq\label{sympl}
\Omega={2\over 1-d} ds\wedge dt^n +\sum_{\alpha=2}^n (\mu_\alpha+{1\over 2})^{-1} dt^{\alpha'2''}\wedge
dt^{(n-\alpha+1)' 1''}.
\eeq
Particularly, the $n$-dimensional planes
\beq\label{lagr}
L_s := \left\{ s={\rm fixed}, ~~t^{\alpha'2''}=0, \alpha=2, \dots, n,
~~t^{\beta'1''}={\rm arbitrary}, ~~\beta=1, \dots, n\right\}
\eeq
are Lagrangian submanifolds in $M\otimes QH^*(CP^1)$.

We will now construct another system $(X^1, \dots, X^n, Y_1, \dots, Y_n)$ of
canonical coordinates using odd periods on $M\otimes QH^*(CP^1)$.
We will choose an appropriate polarization in the space of odd periods on
(\ref{tensor})
and we will compute the generating function $S=S(X,s)$ of the family of
Lagrangian submanifolds
\beq\label{action}
L_s =\left\{ Y_a ={\partial S(X,s)\over \partial X^a}, ~~a=1, \dots, n\right\} .
\eeq
The generating function will be found as an expansion near the point of
``classical limit'' $s=-\infty$. For the particular cases where $M$ is one of
the polynomial Frobenius manifolds of the $ADE$ type (see Section 2 above) we
will identify the generating function with the Seiberg - Witten prepotential of
the four-dimensional supersymmetric Yang - Mills with one of the $ADE$ gauge
groups resp.

Let us spell out the differential equations for the components of odd period
mapping on the tensor product (\ref{tensor}). For any odd period $\varpi$ 
we introduce row vectors $p=(p_1, \dots, p_n)$, $q=(q_1, \dots, q_n)$
$$
p_\alpha={\partial \varpi\over \partial t^{\alpha' 1''}}\equiv {\partial \varpi\over
\partial t^\alpha}, ~~q_\alpha={\partial \varpi\over \partial t^{\alpha' 2''}},
~~q_1 \equiv {\partial \varpi\over \partial s}.
$$

\begin{lemma} At the points of the coordinate cross (\ref{cross}) on $M\otimes
QH^*(CP^1)$ the differential equations (\ref{g-m-l-mu}) with $\nu=1/2$ read
\eqa\label{cross-eq}
&&
\partial_\alpha p\cdot {\mathcal U} + 2 \partial_\alpha q = p(\mu-{1\over 2})
C_\alpha
\nn\\
&&
\partial_\alpha q\cdot {\mathcal U} + 2 Q \partial_\alpha p = q(\mu+{1\over
2})C_\alpha
\eeqa
\beq\label{cross-eq1}
Q\partial_Q (p,q) \left(\begin{matrix} {\mathcal U} & 2\, Q\cr 2 & {\mathcal U}\cr
\end{matrix}\right) =
(p,q) \left( \begin{matrix}0 & Q \, (\mu -{1\over 2})\cr \mu+{1\over 2} & 0\cr
\end{matrix}\right)
\eeq
where $Q=e^s$.
\end{lemma}

Proof is given by straightforward computation using (\ref{g-m-l-mu}),
(\ref{init}), (\ref{mu-tens}), and (\ref{cp1}).

Particularly, near $Q=0$ one can rewrite the equation (\ref{cross-eq1}), for $t\in M\setminus
\Sigma$, in the form
\beq\label{q-lim}
\partial_Q (p,q) ={1\over Q} (p,q) A_0 + O(1)
\eeq
where
$$
A_0 =\left(\begin{matrix} 0 & 0\cr (\mu+{1\over 2}){\mathcal U}^{-1} &
0\cr\end{matrix}\right).
$$
Therefore $Q=0$ is a regular singularity of (\ref{cross-eq1}). The Jordan normal form of the
matrix $A_0$ consists of $n$ $2\times 2$ nilpotent Jordan blocks. So the system
(\ref{cross-eq1}) admits $n$ independent solutions analytic at $Q=0$. They can be completed to
produce a basis by
adding $n$ solutions behaving like $\log Q$ at $Q\to 0$ . We will now explain
how to choose this basis in order to obtain solutions to the full system
(\ref{cross-eq}) -
(\ref{cross-eq1}).

\begin{theorem} Let $\left( x^1(t), \dots, x^n(t)\right)$ be a system
of independent flat coordinates of the intersection form on $M$ defined locally
on $M\setminus \Sigma$. Denote
$$
G^{ab} =(dx^a, dx^b), ~~(G_{ab})=(G^{ab})^{-1}.
$$
Then there exists a basis of odd periods on $M\otimes QH^*(CP^1)$ can be represented in the form
$$
(X^1(t,Q), \dots, X^n(t,Q), Y_1(t,Q), \dots, Y_n(t,Q))
$$
where, at the points of the coordinate cross (\ref{cross}),
\beq\label{lim-x}
X^a =x^a(t) +O(Q), ~~Q\to 0, ~~a=1, \dots, n
\eeq
\beq\label{lim-y}
Y_a =x_a(t)\log Q - 2{\partial F_*(x)\over \partial x^a} + o(1), ~~x_a(t)
=G_{ab}x^b(t), ~a=1, \dots, n.
\eeq
Here $F_*(x)$ is the potential of the dual Frobenius manifold $M_*$. The
functions $X^a$, $Y_b$ are canonically conjugated w.r.t. the symplectic
structure (\ref{sympl}):
\beq\label{canon}
\left\{ X^a, Y_b\right\} =\delta^a_b, ~~\left\{ X^a, X^b\right\} = \left\{ Y_a,
Y_b\right\} =0.
\eeq
The coordinates $X^a$ are determined uniquely, the coordinates $Y_a$ are
determined with the ambiguity that can be absorbed by a redefinition of the 
dual potential
$$
F_*(x) \mapsto F_*(x) + {\rm quadratic}.
$$
\end{theorem}

\pf Left eigenvectors of the matrix $A_0$ are of the form
$$
q=0, ~p={\rm arbitrary}.
$$
Hence for an arbitrary $p^0 =(p^0_1, \dots, p^0_n)$ the system (\ref{cross-eq1}) admits a
solution of the form
$$
(p,q) =(p^0,0)+O(Q), ~~Q\to 0.
$$
The dependence of this solution on $t$ is to be determined from (\ref{cross-eq}).
Particularly, for $p^0$ one obtains the equations coinciding with the system of
differential equations (\ref{g-m}) for the gradients of the flat coordinates of the
intersection form on $M$. This gives the solutions (\ref{lim-x}). The coefficients of the
$Q$-expansion of the solutions are uniquely determined from the system
(\ref{cross-eq}),
(\ref{cross-eq1}):
\eqa\label{instanton}
&&
p_\alpha^{(a)}:= \partial_\alpha X^a =\partial_\alpha x^a(t)
+\partial_\alpha\sum_{k=1}^\infty {Q^k\over (k!)^2} \partial_1^{2\,k}x^a(t)
\nn\\
&&
q_\alpha^{(a)} := \partial_{\alpha' 2''} X^a =\partial_\alpha\sum_{k=1}^\infty {k\, Q^k\over
(k!)^2} \partial_1^{2k-1} x^a(t).
\eeqa

Let us find the second half of the solutions to (\ref{cross-eq}),
(\ref{cross-eq1}). These must have the
form
\beq\label{q-exp}
(p,q) =\left( p^0 \log Q + r^0, q^0\right) + o(1), ~~Q\to 0.
\eeq
Substituting in (\ref{cross-eq1}) we find
\beq\label{q-exp1}
q^0=p^0 {\mathcal U}(\mu+{1\over 2})^{-1}.
\eeq
Now we plug the expansion (\ref{q-exp}), (\ref{q-exp1}) into differential
equations (\ref{cross-eq}). We find, as
before, that $p^0 =(p^0_1, \dots, p^0_n)$ depends on $t$ as
\beq\label{p0}
p^0_\alpha =\partial_\alpha x(t)
\eeq
for some flat function $x(t)$ of the intersection form on $M$. For the dependence
of $r^0$ on $t$ we obtain a system
$$
\partial_\alpha r^0 \cdot {\mathcal U}(t) = r^0 (\mu-{1\over 2}) C_\alpha -2 p^0(t)
C_\alpha(t).
$$
Let us look for a solution of the linear inhomogeneous system in the form
$$
r^0 =A_i(t) \xi^{i}(t), ~~\xi^i(t) =(\partial_1 x^i(t), \dots,
\partial_nx^i(t))
$$
where the coefficients $A_i(t)$ are to be determined. This gives
$$
\partial_\alpha A_i \, \xi^i \, {\mathcal U} =-2 p^0 C_\alpha.
$$
Multiplying both sides by $\eta^{-1}(\xi^c)^T$ we obtain
$$
\partial_\alpha A_i G^{ic} =-2 p^0_\beta c^{\beta\gamma}_\alpha {\partial
x^c\over \partial t^\gamma}.
$$
Applying chain rule we rewrite the last equation in the form
$$
{\partial A_i\over \partial x^b} G^{ic} =-2 p^0_\beta {\partial t^\alpha\over
\partial x^b} c^{\beta\gamma}_\alpha {\partial x^c\over \partial t^\gamma}.
$$
Choosing in (\ref{p0}) $x=x_a(t)$ we obtain
\eqa\label{da}
&&
{\partial A_i\over \partial x^b} 
=-2 G_{ic} G_{am} {\partial x^m\over \partial
t^\beta} {\partial t^\alpha\over \partial x^b} c^{\beta\gamma}_\alpha {\partial
x^c\over \partial t^\gamma}
\nn\\
&&
=-2 {\partial^3 F_*(x)\over \partial x^i \partial x^a \partial x^b}.
\eeqa
So 
$$
A_i =-2 {\partial^2 F_*(x) \over \partial x^i \partial x^a}
$$
and
\eqa\label{r0}
&&
r^0_\alpha =- 2 {\partial^2 F_*(x)\over \partial x^a \partial x^i} {\partial
x^i\over \partial t^\alpha} 
\nn\\
&&
 = -2 {\partial\over \partial t^\alpha} {\partial F_*(x)\over \partial
x^a}.
\eeqa
This gives the solution $Y_a$ of the form (\ref{lim-y})). This solution is determined
uniquely up to adding of a linear combination of the solutions analytic at
$Q=0$.

Let us compute the Poisson brackets of the functions $X$, $Y$. Due to 
(\ref{pb-tens}) we
have, for two arbitary functions $\varpi_1$, $\varpi_2$ on $M\otimes QH^*(CP^1)$
$$
\left\{ \varpi_1, \varpi_2\right\} =<p_1, q_2(\mu+{1\over 2})> + 
<q_1, p_2(\mu-{1\over
2})>.
$$
Here 
$$
(p_i)_\alpha ={\partial \varpi_i\over \partial t^{\alpha'1''}}, ~~(q_i)_\alpha
={\partial \varpi_i\over \partial t^{\alpha'2''}}, ~~\alpha=1, \dots, n, ~~i=1, \,
2.
$$
Due to independence of the brackets of the functions $X^a$ and $Y_b$ on the
point of the Frobenius manifold (\ref{tensor}), we obtain
$$
\{ X^a, X^b\}_{Q=0}=0, ~~
\left\{ X^a, Y_b\right\} =<\xi^{a}, G_{bc}\xi^c{\mathcal U}> =\delta^a_b.
$$
To prove that $\{Y_a, Y_b\}=0$ we will first verify that $\{ \tilde Y_a, \tilde
Y_b\}=0$
where 
$$\tilde Y_a = Y_a+2{\pal F_*(x)\over \pal x^a} = x_a\log Q +o(1).
$$ 
We have
$$
\{\tilde Y_a ,\tilde Y_b\} =\log Q \, G_{ac}G_{bd} 
\left[ <\xi^c ,\xi^d {\mathcal U}> +
<\xi^c{\mathcal U}(\mu+{1\over 2})^{-1}, \xi^d(\mu-{1\over 2})> \right] +o(1).
$$
Using symmetry of ${\mathcal U}$ and skew-symmetry of $\mu$ we obtain that 
the coefficient of $\log Q$ vanishes. Finally we observe that the gradient shift
$$
Y_a =\tilde Y_a -2 {\partial F_*\over \partial x^a}
$$
preserves the Poisson brackets. 
The Theorem is proved.\epf 

We will now describe the generating function (\ref{action}) of the family of Lagrangian
manifolds $L_s$ w.r.t. the canonical coordinates $X^a$, $Y_b$.

\begin{theorem} The generating function $S=S(X,Q)$, $Q=e^s$, of the family
of Lagrangian manifolds $L_s$ admits an expansion
\beq\label{lim-s}
S={1\over 2} G_{ab} X^a X^b \,\log Q-2\, F_*(X)
+\sum_{k\geq 1} S_k(X) Q^k, ~~Q\to 0.
\eeq
The coefficients $S_k(X)$ can be uniquely determined from the Hamilton - Jacobi
equation
\beq\label{ham-jac}
{\partial S\over \partial s} =H
\eeq
with the Hamiltonian
\beq\label{ham}
H={2\over 1-d} t_1 ={1\over 2} G_{ab}x^a x^b
\eeq
and from the expansions
\beq\label{ser-x}
X^a =x^a + \sum_{k=1}^\infty {Q^k\over (k!)^2} \partial_1^{2k} x^a, ~~a=1,
\dots, n.
\eeq
\end{theorem}

\pf The shifts along $s$ are given by a Hamiltonian flow with the linear
Hamiltonian coinciding with the flat coordinate on (\ref{tensor}) conjugate to 
$s$. From
(\ref{sympl})
we find (\ref{ham}) (the expression of the Hamiltonian via the flat coordinates of the
intersection form follows from Lemma \ref{emma}). Substituting the functions
$x^a=x^a(X,Q)$ inverse to (\ref{ser-x}) one obtains an expansion
$$
H={1\over 2} G_{ab}X^a X^b +\sum_{k\geq 1}H_k(X) Q^k
$$
of the Hamiltonian evaluated on the Lagrangian manifold $L_s$. 

The generating function of a family of Lagrangian manifolds $L_s$ moving along
trajectories of a Hamiltonian flow is known to satisfy the Hamilton - Jacobi
equation (\ref{ham-jac}) where in the r.h.s. stands the Hamiltonian of the flow evaluated on
the Lagrangian manifold (see, e.g., \cite{dfn}). This gives
$$
S={1\over 2} G_{ab}X^a X^b \, \log Q +\sum_{k\geq 1} H_k(X){Q^k\over k} +S_0(X).
$$
The integration constant $S_0(X)$ can be determined from the expansion
(\ref{lim-y}) of
$Y_a=\partial S/\partial X^a$. The Theorem is proved.\epf

We give now an integral representation of the odd periods $\varpi(t,Q)$ on
$M\otimes QH^*(CP^1)$ via the even periods on $M$.

\begin{prop} Let $x(t)$ be a period on $M\setminus \Sigma$. Then
\beq\label{elliptic}
\varpi(t,Q) ={1\over 2\pi\,i} \oint {x(t^1-\lambda, t^2, \dots, t^n)\over
\sqrt{\lambda^2-Q}} d\lambda
\eeq
is an odd period on $M\otimes QH^*(CP^1)$ evaluated on the coordinate cross
(\ref{cross}).
\end{prop}

\pf Using (\ref{g-m-l-mu}), (\ref{g-m-l-lam}) and integration by parts we can easily verify that the
functions
$$
p_\alpha = {1\over 2\pi\,i} \oint {\partial_\alpha x(t^1-\lambda, t^2, \dots, t^n)\over
\sqrt{\lambda^2-Q}} d\lambda,
~~q_\alpha ={1\over 2\pi\,i} \oint {\lambda \partial_\alpha x(t^1-\lambda, t^2, \dots, t^n)\over
\sqrt{\lambda^2-Q}} d\lambda
$$
satisfy the system (\ref{cross-eq}), (\ref{cross-eq1}). The Proposition 
is proved.\epf

{\bf Example 1.} In the simple case of one-dimensional Frobenius manifold $M$
one obtains odd period mapping $(t,s)\mapsto (X,Y)$ on $QH^*(CP^1)$. The
expansion (\ref{lim-s}) reads
$$
X=x-\sum_{k\geq 1} {Q^k\over x^{4k-1}} {2^{2k}(4k-3)!!\over (k!)^2}.
$$
The potential $F_*$ is equal to ${1\over 2} x^2 \log x^2$. We obtain therefore
the expansion of the generating function
$$
S={1\over 2} x^2 \log Q - x^2 \log x^2 + ...
$$
The integral (\ref{elliptic}) expresses the periods $X$, $Y$ via complete elliptic integrals
$$
X={2\over \pi i}\int_{-\sqrt Q}^{\sqrt Q} {\sqrt{t-\lambda}\over \sqrt{\lambda^2
-Q}}d\lambda, ~~
Y={2\over \pi i}\int_{-\sqrt Q}^{t} {\sqrt{t-\lambda}\over \sqrt{\lambda^2
-Q}}d\lambda.
$$

{\bf Example 2.} For the polynomial Frobenius manifold $M={\mathbb C}^n/W$ on the
orbit space of a finite Coxeter group (see the end of Section 1) the derivatives
$\partial_1 x^1$, \dots, $\partial_1 x^n$ as functions of $x^1$, \dots, $x^n$
can be found in a pure algebraic way from the linear system
$$
\sum_{a=1}^n \partial_1 x^a {\partial y^i\over \partial x^a} =\delta^i_1, ~~i=1,
\dots, n.
$$
Iterating this procedure we obtain derivatives $\partial^{2k}_1x^a$ for any
$k\geq 1$. This gives an algebraic algorithm of computing of the terms of the
expansion (\ref{lim-s})
\beq\label{cox-s}
S(X) ={1\over 2} G_{ab}X^aX^b\log Q-{1\over 2} \sum_{\alpha\in \Delta_+}
(\alpha,X)^2 \log (\alpha,X)^2
+\sum_{k\geq 1} S_k(X) Q^k.
\eeq

Our observation is that, for the case of $W=$ one of the Weyl groups of the
$ADE$ type, the expansion (\ref{cox-s}) coincides, after redefining 
$$
Q=\Lambda^{4\over 1-d} =\Lambda^{2 h}
$$
$h=$ Coxeter number of $W$, and changing the notations
$$
X\to a, ~~Y\to a_D, ~~S\to {\bf F}
$$
with the instanton expansion of the Seiberg - Witten prepotential of the
4-dimensional $N=2$ supersymmetric Yang - Mills (see in \cite{sw1, sw2, hoker, mmm}). The shortest way
to see this is to eliminate $q$ from the systems (\ref{cross-eq}), (\ref{cross-eq1}). We will do it
assuming semisimplicity of the Frobenius manifold under an additional assumption
$d\neq 3$.

\begin{prop}
The odd periods $\varpi(t,Q)$ evaluated on the coordinate cross (\ref{cross}) 
satisfy the following system of differential equations
\eqa\label{varpi1}
&&
Lie_E \varpi + 2 \partial_s\varpi ={1-d\over 2}\varpi
\\
&&
\partial_s^2\varpi = Q \partial_1^2\varpi\label{varpi2}
\\
&&
\partial_\alpha\partial_\beta\varpi = c_{\alpha\beta}^\gamma(t)
\partial_1\partial_\gamma\varpi.\label{varpi3}
\eeqa
\end{prop}

\pf The equation (\ref{varpi1}) is just a manifestation of the general quasihomogeneity
of the odd periods (see (\ref{quasi-p}) for $\nu={1\over 2}$). To derive 
(\ref{varpi2}) we consider the first component
of the second equation in the system (\ref{cross-eq1}). It can be rewritten, 
due to ${\mathcal U}_1^\alpha=E^\alpha$ in
the form
$$
\partial_\alpha g ={3-d\over 2} q_\alpha
$$
where we denoted
$$
g:= E^\epsilon q_\epsilon +2Q\,p_1.
$$
Substituting this expression into the first equation we obtain, for $p_\epsilon
=\partial_\epsilon\varpi$
$$
\partial_\alpha\partial_\epsilon\varpi \,{\mathcal U}^\epsilon_\beta
-p_\sigma (\mu-{1\over 2})_\epsilon^\sigma c_{\alpha\beta}^\epsilon =-{4\over
3-d} \partial_\alpha\partial_\beta g.
$$
Hence the operator ${\mathcal U}$ of multiplication by $E$ is symmetric w.r.t. the
bilinear form with the matrix 
\beq\label{var-bil}
(\partial_\alpha\partial_\beta \varpi).
\eeq
Any power of ${\mathcal U}$ will still be symmetric w.r.t. the same bilinear form.
On a semisimple Frobenius manifold powers of $E$ span the whole algebra on $T\,
M$.
So, the bilinear form (\ref{var-bil}) must be an invariant form on the Frobenius algebra.
Thus it must have a representation
$$
\partial_\alpha\partial_\beta \varpi =c_{\alpha\beta}^\gamma f_\epsilon
$$
for some covector $f_\epsilon$. Obviously
$$
f_\alpha=\partial_\alpha\partial_1\varpi .
$$
This gives equation (\ref{varpi2}).

To derive (\ref{varpi3}) we rewrite the first component of the equation
$$
2Q\partial_s p + \partial_s q\, {\mathcal U} =Q\, p\, (\mu-{1\over 2})
$$
(see (\ref{cross-eq1})) in the form
$$
2Q\partial_s p_1 + \partial_s(E^\epsilon q_\epsilon) =-Q{1+d\over 2} p_1.
$$
From the definition of the function $g$ it follows that $E^\epsilon q_\epsilon = g -2Q\, p_1$. So
$$
\partial_s g ={3-d\over 2} Q\, p_1.
$$
Differentiating this equation along $t^1$ and using that
$$
\partial_1 g ={3-d\over 2} \partial_s\varpi
$$
we obtain (\ref{varpi3}). Lemma is proved.

For the case of polynomial Frobenius manifolds on the orbit spaces of the $ADE$
Weyl groups the structure constants
$c_{\alpha\beta}^\gamma(t)$ are the same as in (\ref{eguchi}) for the corresponding 
$ADE$ singularity. As it was shown in \cite{ito} the latter equations coincide with
the Picard - Fuchs equations for the periods of Seiberg - Witten differential
on the spectral curve (the latter does not appear in our formalism).
Observe that our Hamilton - Jacobi equation (\ref{ham-jac}) coincides with the
renormalization group equation introduced in \cite{matone, hoker}.

Taking into account that Frobenius manifolds $M$ arise in the setting of 2D
topological field theory, it would be interesting to figure out if there is a
physical motivation for the tensor product $M\otimes QH^*(CP^1)$ to carry 
all the structures of the quantum moduli space (perhaps, all but positivity
of the kinetic energy $\pal^2 S/ \pal X^i \pal X^j$) of a ${\mathcal N}=2$
supersymmetric 4D theory.

\end{document}